\documentclass[11pt]{article}

\usepackage{amssymb,amsthm,amsmath,hyperref,txfonts}

\usepackage{color}

\usepackage{graphicx,float}

\newcommand \nc{\newcommand}
\newtheorem{theorem}{Theorem}[section]
\newtheorem{lemma}[theorem]{Lemma}
\newtheorem{proposition}[theorem]{Proposition}
\newtheorem{corollary}[theorem]{Corollary}
\newtheorem{definition}[theorem]{Definition}

\newtheorem{remark}[theorem]{Remark}

\nc{\ba}{\begin{array}}\nc{\ea}{\end{array}}
\nc{\be}{\begin{eqnarray}}\nc{\ee}{\end{eqnarray}}
\nc{\beq}{\begin{equation}}\nc{\eeq}{\end{equation}}
\nc{\bex}{\begin{eqnarray*}}\nc{\eex}{\end{eqnarray*}}
\nc{\btm}{\begin{theorem}} \nc{\etm}{\end{theorem}}
\nc{\blm}{\begin{lemma}} \nc{\elm}{\end{lemma}}

\newcommand{\sanhao}{\fontsize{16.9pt}{\baselineskip}\selectfont}  

\nc{\R}{\mathbb{R}} \nc{\va}{\varepsilon} \nc{\ls}{\limits}

\def\pf{\noindent{\bf Proof.\quad}}\def\endpf{\hfill$\Box$}

\allowdisplaybreaks

\date{}

\begin{document}
\title{\sanhao On global solutions to a viscous compressible two-fluid model with unconstrained transition to single-phase flow in three dimensions}
\author{Huanyao Wen\footnote{School of Mathematics, South China University of Technology, Guangzhou 510641, China. Corresponding author, e-mail: mahywen@scut.edu.cn.}}
\maketitle

\begin{abstract} We consider the Dirichlet problem for a compressible two-fluid model in multi-dimensions. It consists of the continuity equations for each fluid and the momentum equations for the mixture. This model can be derived from the compressible two-fluid model with equal velocities \cite{BDGG} and from a scaling limit of the Vlasov-Fokker-Planck/compressible Navier-Stokes system \cite{MV3} (see also the compressible Oldroyd-B model with stress diffusion \cite{BLS}). Another interesting connection is that it is formally the equations of compressible magnetohydrodynamic (MHD) flows without resistivity in two dimensions under the action of vertical magnetic field (\cite{Li-Sun}).  Under
weak assumptions on the initial data which can be discontinuous, unbounded and large as well as
involve transition to pure single-phase points or regions, we show existence of global
 weak solutions with finite energy. The essential novelty of this work, compared with
previous works on the same model, is that
transition to each single-phase flow is allowed without any constraints between adiabatic constants or two densities. It means that one of the phases
can vanish in a point while the other can persist. The lack of enough regularity for each densities brings up essential difficulties in the two-component pressure compared with the single-phase model, i.e., compressible Navier-Stokes equations. The key points to achieve the main result rely on the variables reduction technique for the pressure function, domain separation, and some new estimates. As a byproduct, we obtain the existence of global weak solutions to the compressible MHD system without resistivity in two dimensions under the action of non-negatively vertical magnetic field, which represents a step forward to the study of the global large solution to the compressible MHD system without resistivity.
\end{abstract}

{\bf\noindent keywords.} Compressible two-fluid model; global weak
solution; large initial data; transition to each single-phase flow.


\tableofcontents

\section{Introduction}
\setcounter{equation}{0} \setcounter{theorem}{0}
\subsection{Background and motivation}
Multi-phase fluid models have various applications in different areas, such as the petroleum
industry, nuclear, chemical-process, and cryogenics\cite{Br,BDGG,E,EK,I,Z}. They also are quite relevant for the studies of some models like cancer cell migration model \cite{EvjeCES2017,EvjeSIAM}, MHD system \cite{Li-Sun}, and compressible two-fluid Oldroyd-B model with stress diffusion \cite{BLS}.
In this paper, we consider the Dirichlet problem for a viscous compressible two-fluid model with one velocity and a
pressure of two components in three spatial dimensions, i.e.,
    \begin{equation}\label{equation}
    \left\{
    \begin{array}{l}
        n_t+\mathrm{div}(n u)=0,\\
        \rho_t+\mathrm{div}(\rho u)=0,\\
 \big[(\rho+n) u\big]_t+\mathrm{div}\big[(\rho+n) u\otimes u\big]
 +\nabla P(n,\rho)=
            \mu \Delta u + (\mu+\lambda)\nabla \mathrm{div}u \quad \mathrm{on} \   \Omega\times (0, \infty),
    \end{array}
    \right.
    \end{equation}
with  the  initial-boundary conditions
    \begin{equation}\label{initial-boundary}\begin{cases}
      n(x,0)=n_0(x),\ \rho(x,0)=\rho_0(x),\ (\rho+n) u(x,0)=M_0(x)  \ \ \mathrm{for} \ \  x\in\overline{\Omega},
    \\[2mm]
     u|_{\partial\Omega}=0\quad \mathrm{for}\ t\ge0,
    \end{cases}\end{equation}
     where $\Omega\subset\R^3$ is a bounded domain, and $\rho$ and $n$, $u$, and $P$ represent the densities of two fluids, the velocity of the fluids, and the pressure, respectively.
   $\mu$ and $\lambda$ are the viscosity coefficients satisfying
  $\mu> 0,\;\;2\mu+3\lambda\geq 0.$ Here we assume that $\mu$ and $\lambda$ are constants.
  The pressure we study is given by \be\label{pressure law1}
P(n,\rho)=n^\Gamma+\rho^\gamma, \ee or by \be\label{pressure law2}
\begin{cases}P=A_+(\rho_+)^{\gamma}=A_-(\rho_-)^{\Gamma},\\
\rho \rho_-+n\rho_+=\rho_+\rho_-,
 \end{cases}\ee for constants $A_+,A_->0$ and $\gamma, \Gamma>1$, where $\rho=\alpha\rho_+$, $n=(1-\alpha)\rho_-$, and
 $\alpha=\alpha(x,t)\in[0,1]$ denotes the volume fraction of the fluid + in the mixture. For (\ref{pressure law2}), one can use the implicit function theorem to define $\rho_+=\rho_+(n,\rho)$ and  $\rho_-=\rho_-(n,\rho)$ which represent the densities of the fluids $+$ and $-$, respectively  (please see \cite{BDGG,Bresch-Mucha,Novotny} for more details). Note that (\ref{pressure law1}) is motivated by a limiting system derived from Vlasov-Fokker-Planck/compressible Navier-Stokes system \cite{MV3}, and by compressible MHD system for two-dimensional case \cite{Li-Sun}, and by compressible Oldroyd-B model with stress diffusion \cite{BLS}, and that (\ref{pressure law2}) is motivated by the compressible two-fluid model with possibly unequal velocities \cite{BDGG}.

  Our aim is to study the global existence of weak solution to (\ref{equation}) with large initial data in three dimensions. When $\rho\equiv0$ or $n\equiv0$, the system (\ref{equation}) reduces to compressible Navier-Stokes equations for isentropic flow. In this case, some pioneering works on this topics have been achieved. More specifically, Lions \cite{Lions} obtained the first global existence result on weak solution with large initial data in multi-dimensions, where $P=R\rho^\gamma$ for some positive constant $R$ and any given $\gamma\ge\frac{9}{5}$ for three dimensions.
  The constraint for $\gamma$ was relaxed to $\gamma>\frac{3}{2}$ by Feireisl \cite{Feiresil1} and by Feireisl-Novotn\'y-Petzeltov\'a \cite{Feireisl}, and to $\gamma>1$ by Jiang-Zhang \cite{Jiang-Zhang} for spherically symmetric weak solutions. The pressure function in \cite{Lions,Feireisl,Jiang-Zhang} is monotone and convex, which is very essential for the compactness of density. Feireisl \cite{FeireislJDE} extended the result to the case for more general pressure $P(\rho)$ of monotonicity for $\rho\ge \rho_z$. Very recently, Bresch and Jabin \cite{BJ} developed a new method to derive the compactness of the density which does not rely on any monotonicity assumptions on the pressure. It remains largely open\footnote{The problem has been solved by Jiang-Zhang \cite{Jiang-Zhang} for spherically symmetric weak solutions in multi-dimensions.} whether the above results for three dimensions can be extended to the more physical case that the adiabatic constant $\gamma>1$.

 When the pressure is of two components like in (\ref{equation}), it will become more challenging. Some nice properties of one-component pressure are not available any more due to some cross products like $f_1(\rho)f_2(n)$ or even more implicitly $f_3(n,\rho)f_2(n)$ and $f_3(n,\rho)f_4(\rho)$ for some known scalar functions $f_i$, $i=1,2,3,4$. At the first glance, it seems that more regularity on the densities is required to handle the cross products in the context of passing to the limits. These extra regularity properties are, so far, out of reach for large solutions, and the classical techniques cannot be applied directly on (\ref{equation}).

We will give a brief overview for the relevant results on the model (\ref{equation}).
In fact, the studies of the model have been very active for the past few years.  Some global existence results are obtained, however, mostly subject to the case for the domination conditions\footnote{It means that
$n_0\le\overline{c}_0\rho_0$ or $\rho_0\le \overline{c}_1n_0$ for
some positive constants $\overline{c}_0$ and $\overline{c}_1$, which implies that the two fluids are dominated by one of the fluids.}.
\begin{itemize}
\item For the one-dimensional case, Evje and Karlsen \cite{EK} obtained the global existence result on weak solution with large initial data subject to the domination conditions. The one-dimensional properties of the equations implies that the densities of the fluids are bounded for large initial data. This good property is essential to show strong convergence of the densities in the context of the approximation system. The domination condition was removed later by Evje, the author, and Zhu \cite{EWZ} by introducing a new energy equality, which allows transition
to each single-phase flow. For the global existence of small solutions, please refer for instance to \cite{E,EK2,Yao-Zhu,Yao-Zhu2} and the references therein.

\item For the multi-dimensional case, in particular for three dimensions, some new challenges arise due to the multi-dimensional nonlinearity. The boundedness of the densities can not be derived as the one-dimensional case with large initial data. However, with some smallness assumptions, the boundedness of the density and the derivatives of the other quantities arising in the equations can be derived to handle the cross products conveniently, and we refer the readers to \cite{Guo-Yang-Yao, Hao-Li, Yao-Zhang-Zhu}. In a recent work by  Maltese et al. \cite{MMMNPZ}, the authors considered another interesting model with the pressure of two components which can be transformed to the one with one-component pressure, i.e.,
\be\label{MMMNPZ'equa}
\begin{cases}
\rho_t+\mathrm{div}(\rho u)=0,\\[2mm]
        Z_t+\mathrm{div}(Z u)=0,\\[2mm]
 \big(\rho u\big)_t+\mathrm{div}\big(\rho u\otimes u\big)
 +\nabla Z^\gamma=
            \mu \Delta u + (\mu+\lambda)\nabla \mathrm{div}u.
\end{cases}
\ee Thus it makes the approach for compressible Navier-Stokes equations applicable to prove the global existence of weak solutions to (\ref{MMMNPZ'equa}) with large initial data. After that the authors obtained the equivalence between (\ref{MMMNPZ'equa}) and the original system for $\gamma\ge\frac{9}{5}$. But it is not the case for the two-fluid system.

 Very recently, with large initial data and the domination conditions or alternatively with $\Gamma$ and $\gamma$ close enough, Vasseur, the author, and Yu
\cite{Vasseur-Wen-Yu} obtained the global existence of weak solutions to (\ref{equation}) by decomposing the pressure function and deriving a new compactness theorem for transport equations with possible diffusion, where the pressure is determined by the explicit case
(\ref{pressure law1}) for $$\Gamma>\frac{9}{5}\,\, \mathrm{or}\, \, \gamma>\frac{9}{5}.$$ The result with the domination condition was later extended to the case that both $\Gamma$ and $\gamma$ can touch $\frac{9}{5}$ by Novotn\'y and Pokorn\'y \cite{Novotny} where more general pressure laws covering the cases of
both (\ref{pressure law1}) and (\ref{pressure law2}) were considered.

With large initial data but without any domination conditions in multi-dimensions, the global existence theory for weak solutions only holds for the two-fluid
Stokes equations on the
d-dimensional torus $\mathbb{T}^d$ for $d=2,3$. We refer the readers to the seminal work by Bresch, Mucha, Zatorska \cite{Bresch-Mucha}
where the pressure is given by (\ref{pressure law2}). The proof relies on the Bresch-Jabin's new compactness tools for compressible Navier-Stokes equations and the reformulated system
\be\label{Bresch-reform}
\begin{cases}
 R_t+\mathrm{div}(R u)=0,\\[2mm]
        Q_t+\mathrm{div}(Q u)=0,\\[2mm]
-(\lambda+2\mu)\mathrm{div}u+a^+\Big(Z(R,Q)^{\gamma^+}-\{Z(R,Q)^{\gamma^+}\}\Big)=0,\\[2mm]
\mathrm{rot}\, u=0,\,\, \int_{\mathbb{T}^d}u(x,t)\,dx=0,
\end{cases}
\ee
where $R=\rho=\alpha\rho_+$, $Q=n=(1-\alpha)\rho_-$, $\{Z(R,Q)^{\gamma^+}\}=\Big(\int_{\mathbb{T}^d}Z(R,Q)^{\gamma^+}\,dx\Big)/|\mathbb{T}^d|$, $a^+=A_+$, and $\gamma^+=\gamma$.
Note that \cite{Bresch-Mucha} does not need any domination conditions for $\Gamma,\gamma>1$, although
the nonlinear terms $ \big[(\rho+n) u\big]_t$ and $\mathrm{div}\big[(\rho+n) u\otimes u\big]$ in the momentum equations are ignored so that the momentum equations can be transformed to (\ref{Bresch-reform})$_3$.

\end{itemize}

The case without any domination conditions makes the system (\ref{equation}) more realistic in some physical situations and more ``two fluids" properties from mathematical points of view. In this case, however, it is still open whether the global existence of weak solution exists for possibly large initial data in multi-dimensions.
In this paper, we focus on the Dirichlet problem.

\vskip0.3cm

\subsection{Main result}

Note that for each cases of (\ref{pressure law1}) and (\ref{pressure law2}), the pressure $P(n,\rho)$
 satisfies
\be\label{pressure ineq}
\frac{1}{C_0}(n^{\Gamma}+\rho^{\gamma})\le P(n,\rho)\le
C_0(n^{\Gamma}+\rho^{\gamma}) \ee
 for  some positive constant $C_0$. In fact, (\ref{pressure ineq}) is naturally true for the case
    (\ref{pressure law1}).  For the second case (\ref{pressure law2}), we only
    consider the case of $\gamma\ge \Gamma$, since for the other
    case, it is similar.  More specifically, in view of (\ref{pressure law2})$_1$, we obtain that
    \be\label{lowboundrho-}
    \begin{split}
    \rho_-=&(1-\alpha)\rho_-+\alpha\rho_-=n+\alpha(\frac{A_+}{A_-})^\frac{1}{\Gamma}\rho_+^\frac{\gamma}{\Gamma}
    =n+(\frac{A_+}{A_-})^\frac{1}{\Gamma}\rho\rho_+^{\frac{\gamma}{\Gamma}-1}\\
    \ge&n+(\frac{A_+}{A_-})^\frac{1}{\Gamma}\rho^{\frac{\gamma}{\Gamma}},
    \end{split}
    \ee and that
      \be\label{uppboundrho-}
    \begin{split}
    \rho_-=n+(\frac{A_+}{A_-})^\frac{1}{\Gamma}\rho\rho_+^{\frac{\gamma}{\Gamma}-1}
    =n+(\frac{A_-}{A_+})^{-\frac{1}{\gamma}} \rho{\rho_-}
    ^{1-\frac{\Gamma}{\gamma}}
    \le n+\frac{1}{2}\rho_-+c_0\rho^{\frac{\gamma}{\Gamma}}.
    \end{split}
    \ee
    (\ref{lowboundrho-}) and (\ref{uppboundrho-}) imply (\ref{pressure ineq}).

In addition, for any smooth solution of system \eqref{equation},  the following
energy equalities holds for any time $0\leq t\leq T:$
 \be\label{0-energy-inequality}
\begin{split} &\frac{d}{dt}\int_\Omega\Big[\frac{(\rho+n)|u|^2}{2}+G(\rho,n)\Big]\,dx
+\int_\Omega\Big[\mu|\nabla u|^2
+(\mu+\lambda)|\mathrm{div}u|^2\Big]\,dx= 0,
\end{split}
\ee where\be\label{G} G(\rho,n)=\left\{\begin{array}{l}
\frac{n^\Gamma}{\Gamma-1}+\frac{\rho^\gamma}{\gamma-1}, \ \ \ {\rm if} \ \ P\,\, \mathrm{is}\,\, \mathrm{given}\,\, \mathrm{by}\,\, (\ref{pressure law1}), \\
[3mm] P(n,\rho)\Big(\frac{\alpha
}{\gamma-1}+\frac{1-\alpha}{\Gamma-1}\Big), \ \ \ {\rm if} \ \ P
\,\, \mathrm{is}\,\, \mathrm{given}\,\, \mathrm{by}\,\,
(\ref{pressure law2}).
\end{array}
\right. \ee Here (\ref{G})$_1$ is given in \cite{Vasseur-Wen-Yu}
and (\ref{G})$_2$ follows from (\ref{formenergy}) with
$\epsilon,\delta=0$.

Motivated by (\ref{pressure ineq}) and
(\ref{0-energy-inequality}), in order to make the initial energy
is finite, we set the following conditions on the initial data, i.e.,
\begin{equation}
\label{initial restriction-1}
\begin{split}&
\inf\limits_{x\in \Omega}\rho_0\geq0,\quad \inf\limits_{x\in
\Omega}n_0\geq0,\quad \rho_0\in L^\gamma(\Omega),\quad n_0\in
 L^\Gamma(\Omega),
 \end{split}
 \end{equation}
and
\begin{equation}
\label{restriction-2}
 \frac{M_0}{\sqrt{\rho_0+n_0}}\in L^2(\Omega) \;\text{ where }\,\frac{M_0}{\sqrt{\rho_0+n_0}}=0\;\text{ on }\,\{x\in\Omega|\rho_0(x)+n_0(x)=0\},
\end{equation} where $M_0$ is the initial momentum of the mixture given in (\ref{initial-boundary}).

The definition of weak solution in the energy space is given in
the following sense.
\begin{definition}(Global weak solution)
\label{definition of weak soluton} We call $(\rho, n,
u):\Omega\times (0,\infty)\to\mathbb R_+\times\mathbb
R_+\times\mathbb R^3$ a global weak solution of
(\ref{equation})-(\ref{initial-boundary}) if for any $0<T<+\infty$,
\begin{itemize}
\item $ \rho\in L^\infty\big(0,T;L^\gamma(\Omega)\big),\ n\in
L^\infty\big(0,T;L^\Gamma(\Omega)\big),\ \sqrt{\rho+n} u\in
L^\infty\big(0,T;L^2(\Omega)\big), u\in
L^2\big(0,T;H_0^1(\Omega)\big),$

\item  $(\rho, n, u)\ \mathrm{solves}\ \mathrm{the}\
\mathrm{system}\ (\ref{equation})\ \mathrm{in}\
\mathcal{D}^\prime(Q_T),\ \mathrm{where}\
\mathrm{Q_T}=\Omega\times(0,T),$

\item $\big(\rho, n, (\rho+n) u\big)(x,0)=\big(\rho_0(x),n_0(x),
M_0(x)\big),\quad \mathrm{for\ a.e.}\ x\in \Omega,$


\item $(\ref{equation})_1\ \mathrm{and} (\ref{equation})_2\
\mathrm{hold}\ \mathrm{in}\ \mathcal{D}^\prime\big(\mathbb
R^3\times (0,T)\big)\ \mathrm{provided}\ \rho,n,u\ \mathrm{are}\
\mathrm{prolonged}\ \mathrm{to}\ \mathrm{be}\ \mathrm{zero}\
\mathrm{on}\ \mathbb R^3/\Omega,$

\item  the equation (\ref{equation})$_1$ and (\ref{equation})$_2$
are satisfied in the sense of renormalized solutions, i.e.,\bex
\partial_tb(f)+{\rm{div}}\big(b(f)u\big)+[b^\prime(f)f-b(f)]{\rm{div}}u=0\eex
holds in  $\mathcal{D}^\prime(Q_T)$, for any $b\in C^1(\mathbb R)$
such that $b^\prime(z)\equiv0$ for all $z\in \mathbb R$ large
enough, where $f=\rho,n$.
\end{itemize}
\end{definition}


\bigskip

Now we are in the position to state our main result in the paper.
\begin{theorem} \label{th:1.1} For any given $\Gamma\geq\frac{9}{5}$ and $\gamma\ge\frac{9}{5}$. Assume that $\Omega$ is a
 bounded domain in $\R^3$ of class $C^{2+\nu}$ for some $\nu>0$.
Under the conditions of \eqref{initial
restriction-1}-\eqref{restriction-2}, there exists a global weak
solution $(\rho, n, u)$ to (\ref{equation})-(\ref{initial-boundary}).

\end{theorem}

\medskip

\begin{remark}
In Theorem \ref{th:1.1}, the
global weak solution exists for $\Gamma,\gamma\ge\frac{9}{5}$
without any domination conditions, which implies that
transition to each single-phase flow is allowed. In addition, $\Gamma$ and
$\gamma$ are independent within the interval
$[\frac{9}{5},\infty)$, which indicates that it is not necessary for them to stay close to
each other like
\be\label{0-addition}\max\{\frac{3\gamma}{4},\gamma-1,\frac{3(\gamma+1)}{5}\}<\Gamma<\min\{\frac{4\gamma}{3},\gamma+1,\frac{5\gamma}{3}-1\}\ee
 as in
\cite{Vasseur-Wen-Yu} where $\Gamma,\gamma>\frac{9}{5}$ and the
pressure is given by the explicit pressure (\ref{pressure law1}).

Theorem
\ref{th:1.1} provides the first result on the global solution to
the compressible two-fluid system (\ref{equation}) in multi-dimensions without any
 domination conditions and
smallness assumptions for the pressure given by (\ref{pressure law1}) or
by (\ref{pressure law2}). Note that when $\rho\equiv0$ or $n\equiv0$, Theorem \ref{th:1.1} perfectly matches the result of Lions \cite{Lions} for compressible Navier-Stokes equation in a bounded domain of $\mathbb{R}^3$. Lemma \ref{main 2-le:important} is very essential in the proof, which needs $\rho,n\in L^2\big(0,T;L^2(\Omega)\big)$. As a consequence, we require that $\Gamma+\theta_1\ge2$ and $\gamma+\theta_2\ge2$ in Lemma \ref{3-le:h-inofrho}, which yields $\Gamma,\gamma\ge\frac{9}{5}$. Therefore it remains open whether both $\Gamma$ and $\gamma$ can get close to $\frac{3}{2}$ in three dimensions even for the case with domination conditions.

\end{remark}

\begin{remark}
Note that Theorem \ref{th:1.1} is also true for the two-dimensional case, and that the pressure function (\ref{pressure law1}) for $\Gamma=2$  is corresponding to the compressible MHD system without resistivity in two dimensions under the action of vertical magnetic field \cite{Li-Sun}. Thus as a byproduct, we obtain the existence of global weak solutions to the two-dimensional and non-resistive MHD system with non-negatively vertical magnetic field\footnote{For the compressible MHD system with resistivity, the global existence of weak solutions with large initial data has been achieved by Hu, Wang \cite{Hu-Wang}. However, for the case without resistivity, more essential challenges will arise due to the lack of regularity of the magnetic field.}. For the three-dimensional case with the pressure (\ref{pressure law1}) and $\Gamma=2$, it is motivated by compressible Oldroyd-B model with stress diffusion\cite{BLS}.
\end{remark}

The main ingredients in the proof are stated as
follows. As
mentioned in the previous works
\cite{Bresch-Mucha,Novotny,Vasseur-Wen-Yu}, the main challenges
focus on the pressure of two components which brings out some
cross terms between the two densities. Section \ref{sec5} is the
main ingredient in the proof. In fact, in Section \ref{sec5}, the
main point is to prove that $\overline{P(n,\rho)}=P(n,\rho)$ where
$\overline{P(n,\rho)}$ is the weak limit of the approximate
pressure $P(n_\delta,\rho_\delta)$ as $\delta\rightarrow0^+$. It
suffices to establish the strong convergence of $\rho_\delta$ and
$n_\delta$ as $\delta\rightarrow0^+$. To achieve this, it is crucial to
prove that \be\label{0-1}\begin{cases}
\overline{T_k(\rho)}\,\,\overline{P(n,\rho)}
 \le
\overline{T_k(\rho)P(n,\rho)},\\[4mm]
\overline{T_k(n)}\,\,\overline{P(n,\rho)}
 \le
\overline{T_k(n)P(n,\rho)},
\end{cases}
\ee a.e. on $Q_T$, where $T_k$ is a smooth cut-off function for $k=1,2,\cdot\cdot\cdot$. In
Lions-Feireisl's framework for compressible Navier-Stokes equations, the one-component pressure function with monotonicity and convexity gives rise to $$\overline{T_k(\rho)}\,\,\overline{P(\rho)}
 \le
\overline{T_k(\rho)P(\rho)}.$$ But it is not the case for two-fluid system.

Compared with \cite{Vasseur-Wen-Yu}, the new challenge for the proof in the context of allowing unconstrained transition to single-phase flow is to remove (\ref{0-addition}) and allow the two indexes $\Gamma$, $\gamma$ to touch $\frac{9}{5}$. We state the main differences in the proof as below.

\begin{itemize}
\item First, to justify (\ref{0-1}) without (\ref{0-addition}), we can not use the same decomposition of pressure by Vasseur, the author, and Yu (\cite{Vasseur-Wen-Yu}) in the whole domain $Q_T$ any more, i.e.,\bex
P(n_\delta,\rho_\delta)=A^\Gamma(\rho_\delta+n_\delta)^{\Gamma}+B^\gamma
(\rho_\delta+n_\delta)^{\gamma}+\text{remainder} \eex a.e. on $Q_T$, where $(A,B)=(\frac{n}{\rho+n},\frac{\rho}{\rho+n})$ if
$\rho+n\neq0$, since one can not even guarantee the integrability of $(\rho_\delta+n_\delta)^{\gamma}$ and  $(\rho_\delta+n_\delta)^{\Gamma}$ in the whole domain without (\ref{0-addition}) based on the known estimates of $(\rho_\delta,n_\delta)$ in
$\in L^{\gamma+\theta_2}(Q_T)\times L^{\Gamma+\theta_1}(Q_T)$ for $\theta_2=\theta_2(\gamma)$, $\theta_1=\theta_1(\Gamma)$ (see Lemma
\ref{3-le:h-inofrho}), and we do not even have $P(n,\rho)\le\overline{P(n,\rho)}$ for the pressure (\ref{pressure law2}).
In this work we observe that the weighted
functions $A$ and $B$ are able to cancel some possible
oscillation of $\rho_\delta+n_\delta$. As a matter of fact,
$A\rho_\delta$ and $Bn_\delta$ are bounded in
$L^{\Gamma+\theta_1}(Q_T^\prime)$ and in
$L^{\gamma+\theta_2}(Q_T^\prime)$, respectively, for some domain
$Q^\prime_T\subset Q_T$ where the measure of $Q_T/Q^\prime_T$ is
small enough. This can be achieved by obtaining that $n_\delta-Ad_\delta\rightarrow0$ and $\rho_\delta-Bd_\delta\rightarrow0$ strongly in $L^1(Q_T)$ from (\ref{key bound}) for $s=1$. Thus we are able to justify (\ref{0-1}) a.e. on
$Q^\prime_T$ by means of the decomposition of the pressure and the
cut-off functions on $Q^\prime_T$. Finally, by sending the
measure of $Q_T/Q^\prime_T$ to zero, we get (\ref{0-1}).
See Lemma \ref{3-le:4.6} for more details. For the implicit pressure (\ref{pressure law2}), we introduce a new non-decreasing function $G_{A,B}$ (see (\ref{3-increasing function})), i.e., \bex\begin{split}
G_{A,B}(z):=&P(Az,Bz)-\frac{C_2}{\max\{\Gamma,\gamma\}}\Big[(Az)^{\Gamma}
+ (Bz)^{\gamma}\Big],
\end{split}
\eex to connect the implicit pressure with the explicit convex function $\frac{C_2}{\max\{\Gamma,\gamma\}}\Big[(Az)^{\Gamma}
+ (Bz)^{\gamma}\Big]$. The construction of $G_{A,B}$ is inspired by \cite{FeireislJDE} for the compressible Navier-Stokes equations with non-mono pressure of one component.

\item Second, to allow both $\Gamma$ and $\gamma$ to touch $\frac{9}{5}$, which represents a major step forward for the cases of transition to each single-phase flow and of more general pressure law compared with \cite{Vasseur-Wen-Yu} where $\Gamma,\gamma\in (\frac{9}{5},\infty)$ and the pressure (\ref{pressure law1}) is considered only, it is
important to prove that \bex\begin{split}\|T_k(\rho)-
\overline{T_k(\rho)}+T_k(n)-
\overline{T_k(n)}\|_{L^2(Q_T/Q_{T,k})}\rightarrow 0
\end{split}\eex as $k\rightarrow\infty$, where $Q_{T,k}$ is given
by (\ref{3-10}). On the other hand, it is not difficult to justify
\bex\begin{split}
  \|T_k(\rho)- \overline{T_k(\rho)}+T_k(n)-
\overline{T_k(n)}\|_{L^{1}(Q_T/Q_{T,k})}\rightarrow 0
\end{split}
\eex as $k\rightarrow\infty$. Thus by means of the standard
interpolation inequality, it suffices to get the upper bound of
\be\label{0-3}\begin{split} \|T_k(\rho)-
\overline{T_k(\rho)}+T_k(n)-
\overline{T_k(n)}\|_{L^{\Gamma_{min}+1}(Q_T/Q_{T,k})}\end{split}
\ee uniformly for $k$, where
$\Gamma_{min}+1=\min\{\Gamma,\gamma\}+1>2$. In view of that
$\rho_\delta^\Gamma$ and $n_\delta^\gamma$ might not be bounded in
$L^{p_1}(Q_T)$ uniformly for $\delta$ where $p_1>1$, we derive a
new estimate, i.e., \be\label{new esti}
\lim\limits_{\delta\rightarrow0}\|T_k(Ad_\delta) + T_k(Bd_\delta)
-
T_k(Ad)-T_k(Bd)\|^{\Gamma_{min}+1}_{L^{\Gamma_{min}+1}(Q^\prime_T)}\le
C_k\sigma^\frac{K_{min}-1}{K_{min}}+C\ee where
$d_\delta=\rho_\delta+n_\delta$,
$K_{min}=\min\{\frac{\Gamma+\theta_1}{\Gamma},
\frac{\gamma+\theta_2}{\gamma},2\},$ and
$|Q_T/Q^\prime_{T}|\le\sigma$. Here $C$ is independent of
$\sigma$, $\delta$, and $k$, and $C_k$ is independent of $\sigma$
and $\delta$ but may depend on $k$. With the new estimate (\ref{new esti}),
(\ref{0-3}) can be bounded uniformly for $k$. See Lemma
\ref{3-le5.7} for more details.


\end{itemize}
The rest of the paper is organized as follows. In Section \ref{sec2}, we present some useful lemmas which will be used in the proof of Theorem \ref{th:1.1}. In Section \ref{sec3}, as usual we construct an approximation system with artificial viscosity coefficients in both continuity equations and with artificial pressure in the momentum equations. Then we explore a formal energy estimate due to the more complicated pressure (\ref{pressure law2}) and sketch the proof of the global existence of the solution to the approximation system by virtue of the standard Faedo-Galerkin approach.
In Section \ref{sec4}, we pass the quantities to the limits as the artificial viscosity coefficient goes to zero. With the artificial pressure, the pressure given by (\ref{pressure law1}) or (\ref{pressure law2}) has enough integrability. Then we only need to handle the difficulties arising in the implicit pressure (\ref{pressure law2}) compared with our previous work \cite{Vasseur-Wen-Yu}. In Section \ref{sec5}, we take the limits as the coefficient of artificial pressure, i.e., $\delta$, go to zero. It is the last step for the proof. Some new estimates along with some new ideas are obtained in this section.

\section{Some useful tools}\label{sec2}

\begin{lemma}\label{main 2-le:important}
Let $\nu_K\to 0$ as $K\to+\infty,$ and $\nu_K\geq 0$. If
$\varrho^i_K\geq 0$ for i=1, 2, 3, $\cdot\cdot\cdot$, is a
solution to
\begin{equation}
\label{the parabolic for density}
(\varrho^i_K)_t+\mathrm{div}(\varrho^i_Ku_K)=\nu_K\Delta\varrho^i_K,\;\;\varrho^i_K|_{t=0}=\varrho^i_0,
\;\;\nu_K \frac{\partial\varrho^i_K}{\partial
\nu}|_{\partial\Omega}=0,
\end{equation}
 with $C_0\ge1$ independent of $K$ such that
\begin{itemize}

\item   $
  \|\varrho^i_K\|_{L^2(0,T;L^2(\Omega))}+\|\varrho^i_K\|_{L^\infty(0,T;L^{\gamma^i}(\Omega))}\leq C_0,\;\sqrt{\nu_K}\|\nabla\varrho^i_K\|_{L^2(0,T;L^2(\Omega))}\leq
  C_0.
$
  \item $
  \|u_K\|_{L^2(0,T;H^1_0(\Omega))}\leq C_0.
$ \item for any $K>0$ and any $t> 0$:
\begin{equation}
\label{AAAinitial condition for n2 over
density}\int_{\Omega}\frac{(b^i_K)^2}{d_K}\,dx\leq
\int_{\Omega}\frac{(b^i_0)^2}{d_0}\,dx,
\end{equation}
\end{itemize}
where $b^i_K=\varrho^i_K$, $d_K=\sum\limits_{i=1}^N \varrho^i_K$ for any fixed integer $N\ge2$,
and $\gamma^i>1$.

Then, up to a subsequence, we have
$$\varrho^i_K\to \varrho^i,\;\text{ weakly in } L^2(0,T;L^2(\Omega))\cap L^\infty(0,T;L^{\gamma^i}(\Omega)),$$
$$ u_K\to u \text{ weakly in } L^2(0,T;H^1_0(\Omega)),$$
as $K\rightarrow\infty$, and for any $s\ge1$,
\begin{equation}
\label{key bound}
\lim_{K\to+\infty}\int_0^T\int_{\Omega}d_K|a^i_K-a^i|^s\,dx\,dt=0,
\end{equation}
where $a^i_K=\frac{b^i_K}{d_K}$ if $d_K\neq0$, $a^i=\frac{b^i}{d}$
if $d\neq0$, and $a^i_Kd_K=b^i_K$, $a^id=b^i$ for $i$=1, 2, 3,
$\cdot\cdot\cdot$. Here $(b^i,d)$ is the weak limit of
$(b^i_K,d_K)$ as $K\rightarrow\infty$.
\end{lemma}
\begin{remark}
For $i=1,2$, Lemma \ref{main 2-le:important} can be found in
\cite{Vasseur-Wen-Yu}. It is not difficult to verify the more
general case for $i=1,2,3,\cdot\cdot\cdot$, since (\ref{the
parabolic for density}) is a linear equation.
The compactness conclusion here for the multi-equations with possible diffusion can be applied to study the multi-fluid system introduced in \cite{Novotny} where $P=P(\rho_1,\rho_2,\cdot\cdot\cdot, \rho_N)$. Note that the proof in \cite{Vasseur-Wen-Yu} relies on the DiPerna-Lions renormalized argument for transport equations \cite{DL1, DL}. Thus the $L^2$ bounds of the densities make it possible to use this theory for equations (\ref{the parabolic for density}).
\end{remark}

\begin{lemma}\cite{Vasseur-Wen-Yu}
\label{main lemma}Let $\beta:\R^N\to \R$ be a $C^1$ function with
$|\nabla\beta(X)|\in L^{\infty}(\R^N)$, and $R\in
\left(L^2(0,T;L^{2}(\Omega))\right)^N,$ $u\in
L^2(0,T;H^1_0(\Omega))$ satisfy
\begin{equation}\label{R-equation}
\frac{\partial}{\partial_t}R+\mathrm{div}(u \otimes R)=0,\;\;
R|_{t=0}=R_0(x)
\end{equation} in the distribution sense.
Then we have
\begin{equation}
\label{renormalized} (\beta(R))_t+\mathrm{div}(\beta(R)
u)+[\nabla\beta(R)\cdot R-\beta(R)]\mathrm{div}u=0
\end{equation}
in the distribution sense. Moreover, if $R\in
L^{\infty}(0,T;L^{\gamma}(\Omega))$ for $\gamma>1$, then
$$R\in C([0,T];L^1(\Omega)),$$ and so
$$\int_{\Omega}\beta(R)\,dx(t)=\int_{\Omega}\beta(R_0)\,dx-\int_0^t\int_{\Omega}[\nabla\beta(R)\cdot
R-\beta(R)]\mathrm{div}u\,dx\,dt.$$
\end{lemma}

\begin{remark}
$N$ in Lemma \ref{main lemma} is specified to be 2 in the present paper (see Lemma \ref{3-le:4.6}).
\end{remark}

\begin{lemma}\label{Feireisl-Novotny lemma}( \cite{Feireisl-Novotny}, Theorem 10.19)
Let $I\subset \mathbb{R}$ be an interval, $Q\subset\mathbb{R}^N$ be a
domain, and
$$
(P,G)\in C(I)\times C(I)\,\, \mathrm{be\, a\, couple\, of\,
non-decreasing\, functions}.
$$
Assume that $\varrho_n\in L^1(Q;I)$ is a sequence of functions such that \bex
\begin{cases}
P(\varrho_n)\rightarrow\overline{P(\varrho)},\\[2mm]
G(\varrho_n)\rightarrow\overline{G(\varrho)},\\[2mm]
P(\varrho_n)G(\varrho_n)\rightarrow\overline{P(\varrho)G(\varrho)},
\end{cases}
\eex weakly in $L^1(Q)$. Then \bex
\overline{P(\varrho)}\,\overline{G(\varrho)}\le
\overline{P(\varrho)G(\varrho)},\,\, \mathrm{a.e.\,\, in\,\, Q}.\eex

\end{lemma}

\section{Existence of solutions to an approximate system}\label{sec3}
\setcounter{equation}{0} \setcounter{theorem}{0} In this section,
we construct a sequence of global weak solution $(\rho,n,u)$ to
the following approximation system \eqref{a-equation}-\eqref{a-boundary}.
Motivated by the work of \cite{Feireisl,Vasseur-Wen-Yu}, we
consider the following approximation system
\begin{equation}\label{a-equation}
    \left\{
    \begin{array}{l}
        n_t+\mathrm{div}(n u)=\epsilon\Delta n,\\[2mm]
        \rho_t+\mathrm{div}(\rho u)=\epsilon\Delta \rho,\\[2mm]
 \big[(\rho+n) u\big]_t+\mathrm{div}\big[(\rho+n) u\otimes u\big]
 +\nabla P(n,\rho)+\delta\nabla(\rho+n)^\beta+\epsilon\nabla u\cdot\nabla(\rho+n)\\ =
            \mu \Delta u + (\mu+\lambda)\nabla \mathrm{div}u
    \end{array}
    \right.
    \end{equation}
    on  $\Omega\times (0, \infty)$, with initial and boundary condition
\be\label{a-initial} \big(\rho, n,(\rho+n)
u\big)|_{t=0}=(\rho_{0,\delta},n_{0,\delta},M_{0,\delta})\
\mathrm{on}\ \overline{\Omega}, \ee \be\label{a-boundary}
(\frac{\partial\rho}{\partial \nu},\frac{\partial n}{\partial\nu},
u)|_{\partial\Omega}=0, \ee where $\epsilon,\delta>0$,
$\beta>\max\{4,\Gamma+1,\gamma+1\}$,
$M_{0,\delta}=(\rho_{0,\delta}+n_{0,\delta})u_{0,\delta}$ and
$n_{0,\delta},\rho_{0,\delta}\in C^3(\overline{\Omega})$,
$u_{0,\delta}\in C_0^3(\Omega)$ satisfies
\be\label{appinitial}\begin{cases} 0<\delta\le
\rho_{0,\delta},n_{0,\delta}\le \delta^{-\frac{1}{2\beta}},\quad
 (\frac{\partial n_{0,\delta}}{\partial
\nu}, \frac{\partial
\rho_{0,\delta}}{\partial \nu})|_{\partial\Omega}=0,\\[2mm]
\lim\limits_{\delta\rightarrow0}\left(\|\rho_{0,\delta}-\rho_0\|_{L^\gamma(\Omega)}+\|n_{0,\delta}-n_0\|_{L^\Gamma(\Omega)}\right)=0,\\[2mm]
u_{0\delta}=\frac{\varphi_\delta}{\sqrt{\rho_{0,\delta}+n_{0,\delta}}}\eta_\delta*(\frac{M_0}{\sqrt{\rho_0+n_0}}),
\\[2mm]
\sqrt{\rho_{0,\delta}+n_{0,\delta}}u_{0,\delta}\rightarrow\frac{M_0}{\sqrt{\rho_0+n_0}}\quad
\mathrm{in}\
L^2(\Omega) \ \mathrm{as}\ \delta\rightarrow0,\\[2mm]
M_{0,\delta}\rightarrow M_0\quad \mathrm{in}\ L^1(\Omega) \
\mathrm{as}\ \delta\rightarrow0,
\end{cases} \ee where $\delta\in(0,1)$, $\eta$ is a standard mollifier, $\varphi_\delta\in C_0^\infty(\Omega)$, $0\le\varphi_\delta\le1$ on
$\overline{\Omega}$ and $\varphi_\delta\equiv 1$ on
$\big\{x\in\Omega| \mathrm{dist}(x,\partial\Omega)>\delta\big\}$.


In order to simplify the presentation of the proof, we only
consider the more complicated case of pressure, i.e.,
(\ref{pressure law2}), in the rest of the paper.

\subsection{A formal energy estimate}
The main difference between the approximation system
(\ref{a-equation})-(\ref{a-boundary}) and the one in
\cite{Vasseur-Wen-Yu} by Vasseur, the author, and Yu is
that one of the pressure functions, i.e., (\ref{pressure
law2}), is more complicated. Therefore we will give a formal energy
estimate in this part so that the Galerkin approach could work as in \cite{Vasseur-Wen-Yu}. More specifically, we consider
the pressure given by (\ref{pressure law2}), and suppose that the
solution to (\ref{a-equation})-(\ref{a-boundary}) is smooth
enough.

Define $P_1(\rho_+)=A_+(\rho_+)^\gamma$ and
$P_2(\rho_-)=A_-(\rho_-)^\Gamma$. Since
$P(n,\rho)=P_1(\rho_+)=P_2(\rho_-)$, we decompose the pressure
into two parts, i.e., \be\label{decom of pressure}
P(n,\rho)=\alpha P(n,\rho) + (1-\alpha) P(n,\rho)=\alpha
P_1(\rho_+) + (1-\alpha) P_2(\rho_-), \ee where
$\alpha=\frac{\rho}{\rho_+}$.
 Actually, the idea for
the decomposition (\ref{decom of pressure}) has been used by Evje, the author, Zhu \cite{EWZ} and by Bresch,
Mucha, Zatorska \cite{Bresch-Mucha} to study the one-dimensional case for the
full compressible two-fluid equations with
singular pressure gradient and multi-dimensional case for the
compressible two-fluid Stokes equations, respectively. It is
motivated by the full compressible two-fluid system with unequal
velocities, see \cite{BDGG,BHL}. However, the Laplacian of $\rho$
and $n$ in (\ref{a-equation}) will make the estimates more
complicated.

Multiplying (\ref{a-equation})$_3$ by $u$, integrating by parts
over $\Omega$, and using (\ref{a-equation})$_1$ and
(\ref{a-equation})$_2$, we have \be\label{Form1}\begin{split}
&\frac{d}{dt}\int_\Omega\frac{1}{2}(\rho+n) |u|^2\,dx+
 +\int_\Omega \Big(\mu|\nabla
 u|^2+(\mu+\lambda)|\mathrm{div}u|^2\Big)\,dx\\=&-\delta\int_\Omega u\cdot\nabla(\rho+n)^\beta\,dx-\int_\Omega u\cdot\nabla P(n,\rho)\,dx\\
 =&I_1+I_2.
\end{split}
\ee For $I_1$, we have \be\label{I1}
\begin{split}
I_1=&
\delta\int_\Omega
\frac{\beta}{\beta-1}(\rho+n)^{\beta-1}\nabla\cdot[(\rho+n)u]\,dx
\\=&
-\delta\int_\Omega
\frac{\beta}{\beta-1}(\rho+n)^{\beta-1}(\rho+n)_t\,dx+\delta\epsilon\int_\Omega
\frac{\beta}{\beta-1}(\rho+n)^{\beta-1}\Delta(\rho+n)\,dx\\=&
-\frac{\delta}{\beta-1}\frac{d}{dt}\int_\Omega
(\rho+n)^{\beta}\,dx-\delta\epsilon\int_\Omega
\beta(\rho+n)^{\beta-2}|\nabla(\rho+n)|^2\,dx.
\end{split}
\ee For $I_2$, by virtue of the decomposition (\ref{decom of
pressure}), we have\bex\begin{split} I_2=&-\int_\Omega
\alpha\nabla P_1(\rho_+)\cdot u\,dx-\int_\Omega (1-\alpha)\nabla
P_2(\rho_-)\cdot u\,dx\\=& I_{2,1}+I_{2,2}.
\end{split}
\eex For $I_{2,1}$, we have \be\label{I21}\begin{split}
I_{2,1}=&-\int_\Omega \frac{\gamma A_+}{\gamma-1}(\rho
u)\cdot\nabla \rho_+^{\gamma-1}\,dx\\=& -\int_\Omega \frac{\gamma
A_+}{\gamma-1}\rho_t\rho_+^{\gamma-1}\,dx+\epsilon\int_\Omega
\frac{\gamma
A_+}{\gamma-1}\rho_+^{\gamma-1}\Delta\rho\,dx\\=&-\frac{d}{dt}\int_\Omega
\frac{\gamma A_+}{\gamma-1}\rho\rho_+^{\gamma-1}\,dx+\int_\Omega
A_+\gamma\rho\rho_+^{\gamma-2}(\rho_+)_t\,dx-\epsilon\int_\Omega
A_+\gamma\rho_+^{\gamma-2}\nabla\rho\cdot\nabla\rho_+\,dx\\=&-\frac{d}{dt}\int_\Omega
\frac{\gamma A_+}{\gamma-1}\rho\rho_+^{\gamma-1}\,dx+\int_\Omega
\alpha(A_+\rho_+^{\gamma})_t\,dx-\epsilon\int_\Omega
A_+\gamma\rho_+^{\gamma-2}\nabla\rho\cdot\nabla\rho_+\,dx.
\end{split}
\ee

Similarly, for $I_{2,2}$, we have \be\label{I22}\begin{split}
I_{2,2}=&-\frac{d}{dt}\int_\Omega \frac{\Gamma
A_-}{\Gamma-1}n\rho_-^{\Gamma-1}\,dx+\int_\Omega
(1-\alpha)(A_-\rho_-^{\Gamma})_t\,dx-\epsilon\int_\Omega
A_-\Gamma\rho_-^{\Gamma-2}\nabla n\cdot\nabla\rho_-\,dx.
\end{split}
\ee
(\ref{I21}) and (\ref{I22}) yield that
\be\label{sumI2122}\begin{split}
I_{2,1}+I_{2,2}=&-\frac{d}{dt}\int_\Omega \Big(\frac{\gamma
A_+}{\gamma-1}\alpha\rho_+^{\gamma}+\frac{\Gamma
A_-}{\Gamma-1}(1-\alpha)\rho_-^{\Gamma}\Big)\,dx+\int_\Omega
(A_+\rho_+^{\gamma})_t\,dx\\&-\epsilon\int_\Omega
\Big(A_+\gamma\rho_+^{\gamma-2}\nabla\rho\cdot\nabla\rho_++A_-\Gamma\rho_-^{\Gamma-2}\nabla
n\cdot\nabla\rho_-\Big)\,dx\\=&-\frac{d}{dt}\int_\Omega
A_+\rho_+^{\gamma}\Big(\frac{\gamma
}{\gamma-1}\alpha+\frac{\Gamma}{\Gamma-1}(1-\alpha)-1\Big)\,dx\\&-\epsilon\int_\Omega
\Big(A_+\gamma\rho_+^{\gamma-2}\nabla\rho\cdot\nabla\rho_++A_-\Gamma\rho_-^{\Gamma-2}\nabla
n\cdot\nabla\rho_-\Big)\,dx\\=&-\frac{d}{dt}\int_\Omega
A_+\rho_+^{\gamma}\Big(\frac{\alpha
}{\gamma-1}+\frac{1-\alpha}{\Gamma-1}\Big)\,dx\\&-\epsilon\int_\Omega
\Big(A_+\gamma\rho_+^{\gamma-2}\nabla\rho\cdot\nabla\rho_++A_-\Gamma\rho_-^{\Gamma-2}\nabla
n\cdot\nabla\rho_-\Big)\,dx,
\end{split}
\ee where we have used $A_-\rho_-^{\Gamma}=A_+\rho_+^{\gamma}$. We
still need to analyze the last integral on the right hand side of
(\ref{sumI2122}). More specifically, substituting
$\rho_-=\big(\frac{A_+}{A_-}\big)^\frac{1}{\Gamma}\rho_+^\frac{\gamma}{\Gamma}$
 into (\ref{pressure
law2})$_2$, and differentiating the result with respect to $x$, we
have
\bex\begin{split}
\big(\frac{A_+}{A_-}\big)^\frac{1}{\Gamma}\nabla\rho\rho_+^\frac{\gamma}{\Gamma}+\big(\frac{A_+}{A_-}\big)^\frac{1}{\Gamma}
\frac{\gamma}{\Gamma}\rho\rho_+^{\frac{\gamma}{\Gamma}-1}\nabla\rho_++\nabla
n\rho_++n\nabla\rho_+=\big(\frac{A_+}{A_-}\big)^\frac{1}{\Gamma}(\frac{\gamma}{\Gamma}+1)\rho_+^{\frac{\gamma}{\Gamma}}\nabla\rho_+,
\end{split}
\eex which implies that\bex\begin{split}
\nabla\rho_+=&\Big[\big(\frac{A_+}{A_-}\big)^\frac{1}{\Gamma}\nabla\rho
\rho_+^\frac{\gamma}{\Gamma}+\nabla
n\rho_+\Big]\Big[\big(\frac{A_+}{A_-}\big)^\frac{1}{\Gamma}(\frac{\gamma}{\Gamma}+1)\rho_+^{\frac{\gamma}{\Gamma}}
-\big(\frac{A_+}{A_-}\big)^\frac{1}{\Gamma}\frac{\gamma}{\Gamma}\rho\rho_+^{\frac{\gamma}{\Gamma}-1}-n\Big]^{-1}.
\end{split}
\eex Note that \bex
n=(1-\alpha)\rho_-=(1-\alpha)\big(\frac{A_+}{A_-}\big)^\frac{1}{\Gamma}\rho_+^\frac{\gamma}{\Gamma}.
\eex Hence we have \be\label{narho+}\begin{split}
\nabla\rho_+=&\Big[\big(\frac{A_+}{A_-}\big)^\frac{1}{\Gamma}\nabla\rho
\rho_+^\frac{\gamma}{\Gamma}+\nabla
n\rho_+\Big]\Big[\big(\frac{A_+}{A_-}\big)^\frac{1}{\Gamma}(\frac{\gamma}{\Gamma}+1)\rho_+^{\frac{\gamma}{\Gamma}}
-\big(\frac{A_+}{A_-}\big)^\frac{1}{\Gamma}\frac{\gamma}{\Gamma}\alpha\rho_+^{\frac{\gamma}{\Gamma}}
-(1-\alpha)\big(\frac{A_+}{A_-}\big)^\frac{1}{\Gamma}\rho_+^\frac{\gamma}{\Gamma}\Big]^{-1}
\\=&\frac{\big(\frac{A_+}{A_-}\big)^\frac{1}{\Gamma}\nabla\rho
\rho_+^\frac{\gamma}{\Gamma}+\nabla
n\rho_+}{(\frac{A_+}{A_-}\big)^\frac{1}{\Gamma}\rho_+^{\frac{\gamma}{\Gamma}}}\Big[\big(\frac{\gamma}{\Gamma}+1)
-\frac{\gamma}{\Gamma}\alpha
-(1-\alpha)\Big]^{-1}\\=&\frac{\big(\frac{A_+}{A_-}\big)^\frac{1}{\Gamma}\nabla\rho
\rho_+^\frac{\gamma}{\Gamma}+\nabla
n\rho_+}{(\frac{A_+}{A_-}\big)^\frac{1}{\Gamma}\rho_+^{\frac{\gamma}{\Gamma}}\Big[\frac{\gamma}{\Gamma}(1-\alpha)
 +\alpha\Big]}.
\end{split}
\ee Since \bex
\rho_-=\big(\frac{A_+}{A_-}\big)^\frac{1}{\Gamma}\rho_+^\frac{\gamma}{\Gamma},
\eex we have \be\label{narho-}
\begin{split}
\nabla\rho_-=&\big(\frac{A_+}{A_-}\big)^\frac{1}{\Gamma}\frac{\gamma}{\Gamma}\rho_+^{\frac{\gamma}{\Gamma}-1}\nabla\rho_+
\\=&\frac{\gamma}{\Gamma}\frac{\big(\frac{A_+}{A_-}\big)^\frac{1}{\Gamma}\nabla\rho
\rho_+^{\frac{\gamma}{\Gamma}-1}+\nabla
n}{\frac{\gamma}{\Gamma}(1-\alpha)
 +\alpha},
\end{split}
\ee where we have used (\ref{narho+}).

Now we are in a position to evaluate the last integral on the right
hand side of (\ref{sumI2122}).  In view of (\ref{narho+}) and
(\ref{narho-}), we have \bex\begin{split} &
A_+\gamma\rho_+^{\gamma-2}\nabla\rho\cdot\nabla\rho_++A_-\Gamma\rho_-^{\Gamma-2}\nabla
n\cdot\nabla\rho_-
 \\=&
A_+\gamma\rho_+^{\gamma-2}\frac{\big(\frac{A_+}{A_-}\big)^\frac{1}{\Gamma}\rho_+^\frac{\gamma}{\Gamma}|\nabla\rho|^2
+\rho_+\nabla\rho\cdot\nabla
n}{(\frac{A_+}{A_-}\big)^\frac{1}{\Gamma}\rho_+^{\frac{\gamma}{\Gamma}}\Big[\frac{\gamma}{\Gamma}(1-\alpha)
 +\alpha\Big]} + A_-\gamma\rho_-^{\Gamma-2}\frac{\big(\frac{A_+}{A_-}\big)^\frac{1}{\Gamma}\rho_+^{\frac{\gamma}{\Gamma}-1}\nabla
n\cdot\nabla\rho +|\nabla n|^2}{\frac{\gamma}{\Gamma}(1-\alpha)
 +\alpha}\\=&\frac{\gamma A_+\rho_+^{\gamma}}{\frac{\gamma}{\Gamma}(1-\alpha)
 +\alpha}\Big[
\rho_+^{-2}\frac{\big(\frac{A_+}{A_-}\big)^\frac{1}{\Gamma}\rho_+^\frac{\gamma}{\Gamma}|\nabla\rho|^2
+\rho_+\nabla\rho\cdot\nabla
n}{(\frac{A_+}{A_-}\big)^\frac{1}{\Gamma}\rho_+^{\frac{\gamma}{\Gamma}}}
+
\rho_-^{-2}\Big(\big(\frac{A_+}{A_-}\big)^\frac{1}{\Gamma}\rho_+^{\frac{\gamma}{\Gamma}-1}\nabla
n\cdot\nabla\rho +|\nabla n|^2\Big)\Big]\\=&\frac{\gamma
A_+\rho_+^{\gamma}}{\frac{\gamma}{\Gamma}(1-\alpha)
 +\alpha}\Big[
\rho_+^{-2}|\nabla\rho|^2 +\frac{\rho_+^{-1}\nabla\rho\cdot\nabla
n}{(\frac{A_+}{A_-}\big)^\frac{1}{\Gamma}\rho_+^{\frac{\gamma}{\Gamma}}}
+\big(\frac{A_+}{A_-}\big)^{-\frac{1}{\Gamma}}\rho_+^{-\frac{\gamma}{\Gamma}-1}\nabla
n\cdot\nabla\rho
+(\frac{A_+}{A_-}\big)^{-2\frac{1}{\Gamma}}\rho_+^{-2\frac{\gamma}{\Gamma}}|\nabla
n|^2\Big]
\\=&\frac{\gamma
A_+\rho_+^{\gamma}}{\frac{\gamma}{\Gamma}(1-\alpha)
 +\alpha}\Big|
\rho_+^{-1}\nabla\rho
+(\frac{A_+}{A_-}\big)^{-\frac{1}{\Gamma}}\rho_+^{-\frac{\gamma}{\Gamma}}\nabla
n\Big|^2.
\end{split}
\eex This combined with (\ref{sumI2122}) yields
\be\label{I2}\begin{split}
I_2&=I_{2,1}+I_{2,2}\\&=-\frac{d}{dt}\int_\Omega
A_+\rho_+^{\gamma}\Big(\frac{\alpha
}{\gamma-1}+\frac{1-\alpha}{\Gamma-1}\Big)\,dx-\epsilon\int_\Omega
\frac{\gamma A_+\rho_+^{\gamma}}{\frac{\gamma}{\Gamma}(1-\alpha)
 +\alpha}\Big|
\rho_+^{-1}\nabla\rho
+(\frac{A_+}{A_-}\big)^{-\frac{1}{\Gamma}}\rho_+^{-\frac{\gamma}{\Gamma}}\nabla
n\Big|^2\,dx.
\end{split}
\ee Combining (\ref{Form1}), (\ref{I1}) and (\ref{I2}), we have
 \be\label{formenergy}\begin{split} &\frac{d}{dt}\int_\Omega\Big[\frac{1}{2}(\rho+n)
|u|^2+
\frac{\delta}{\beta-1}(\rho+n)^{\beta}+A_+\rho_+^{\gamma}\Big(\frac{\alpha
}{\gamma-1}+\frac{1-\alpha}{\Gamma-1}\Big)\Big]\,dx+ \int_\Omega
\Big(\mu|\nabla
 u|^2+(\mu+\lambda)|\mathrm{div}u|^2\Big)\,dx\\=&-\epsilon\int_\Omega
\Big[\delta\beta(\rho+n)^{\beta-2}|\nabla(\rho+n)|^2+\frac{\gamma
A_+\rho_+^{\gamma}}{\frac{\gamma}{\Gamma}(1-\alpha)
 +\alpha}\big|
\rho_+^{-1}\nabla\rho
+(\frac{A_+}{A_-}\big)^{-\frac{1}{\Gamma}}\rho_+^{-\frac{\gamma}{\Gamma}}\nabla
n\big|^2\Big]\,dx.
\end{split}
\ee

\subsection{Faedo-Galerkin approach}

In this part, motivated by \cite{Lions,Feireisl} (see also
\cite{Vasseur-Wen-Yu}), we will use Faedo-Galerkin approach to
construct a global weak solution to (\ref{a-equation}),
(\ref{a-initial}) and (\ref{a-boundary}). To begin with, we
consider a sequence of finite dimensional spaces \bex
X_k=[span\{\psi_j\}_{j=1}^k]^3,\quad
k\in\{1,2,3,\cdot\cdot\cdot\}, \eex where
$\{\psi_i\}_{i=1}^\infty$ is the set of the eigenfunctions of the
Laplacian: \bex\begin{cases} -\Delta\psi_i=\lambda_i\psi_i \quad \mathrm{on}\ \Omega,\\[2mm]
\psi_i|_{\partial\Omega}=0.
\end{cases}\eex

For any given $\epsilon,\delta>0$, we shall look for the
approximate solution $u_k\in C([0,T];X_k)$ (for any fixed $T>0$)
given by the following form: \be\label{1-apu}\begin{split}
&\int_\Omega(\rho_k+n_k) u_k(t)\cdot\psi\,dx-\int_\Omega
m_{0,\delta}\cdot\psi\,dx=\int_0^t\int_\Omega\big[\mu \Delta u_k +
(\mu+\lambda)\nabla
\mathrm{div}u_k\big]\cdot\psi\,dx\,ds\\&-\int_0^t\int_\Omega\Big[\mathrm{div}\big[(\rho_k+n_k)
u_k\otimes u_k\big]
 +\nabla P(n_k, \rho_k)+\delta\nabla(\rho_k+n_k)^\beta+\epsilon\nabla
 u_k\cdot\nabla(\rho_k+n_k)\Big]\cdot\psi\,dx\,ds
\end{split}
\ee for $t\in[0, T]$ and $\psi\in X_k$, where $\rho_k=\rho_k(u_k)$
and $n_k=n_k(u_k)$ satisfy \be\label{1-apnrho}\begin{cases}
\partial_tn_k+\mathrm{div}(n_k u_k)=\epsilon\Delta n_k,\\[2mm]
        \partial_t\rho_k+\mathrm{div}(\rho_k u_k)=\epsilon\Delta \rho_k,\\[2mm]
        n_k|_{t=0}=n_{0,\delta},\quad
        \rho_k|_{t=0}=\rho_{0,\delta},\\[2mm]
(\frac{\partial\rho_k}{\partial \nu},\frac{\partial
n_k}{\partial\nu})|_{\partial\Omega}=0.
\end{cases}
\ee


Due to Lemmas 2.1 and 2.2 in \cite{Feireisl}, the problem
(\ref{1-apu}) can be solved on a short time interval $[0, T_k]$
for $T_k\le T$ by a standard fixed point theorem on the Banach
space $C([0, T_k]; X_k)$. To show that $T_k=T$, as in
\cite{Feireisl} (see also \cite{Vasseur-Wen-Yu}), we only need to
get the energy estimate (\ref{formenergy}) with $(\rho,n,u)$
replaced by $(\rho_k,n_k,u_k)$, which could be done by
differentiating (\ref{1-apu}) with respect to time, taking
$\psi=u_k(t)$ and using (\ref{1-apnrho}). We refer the readers to
\cite{Feireisl} for more details. Thus, we obtain a solution
$(\rho_k, n_k, u_k)$ to \eqref{1-apu}-\eqref{1-apnrho} globally in
time with the following bounds \be\label{1-0}\begin{cases}
0<\frac{1}{c_k}\le \rho_k(x,t),n_k(x,t)\le c_k\ \mathrm{for}\
\mathrm{a.e.} (x,t)\in\Omega\times (0,T), \\[3mm]
\sup\limits_{t\in[0,T]}\|(\rho_{+,k},\rho_k)(t)\|_{L^\gamma(\Omega)}^\gamma\le
C(\rho_0,n_0,M_0),\\[3mm]
\sup\limits_{t\in[0,T]}\|(\rho_{-,k},n_k)(t)\|_{L^\Gamma(\Omega)}^\Gamma\le
C(\rho_0,n_0,M_0),\\[3mm]
\delta\sup\limits_{t\in[0,T]}\|\rho_k(t)+n_k(t)\|_{L^\beta(\Omega)}^\beta\le
C(\rho_0,n_0,M_0),\\[3mm]
\sup\limits_{t\in[0,T]}\|\sqrt{\rho_k+n_k}(t)u_k(t)\|_{L^2(\Omega)}^2\le
C(\rho_0,n_0,M_0),\\[4mm]
\displaystyle\int_0^T\|u_k(t)\|_{H^1_0(\Omega)}^2\,dt\le C(\rho_0,n_0,M_0),\\[4mm]
\displaystyle\epsilon\int_0^T\big(\|\nabla\rho_k(t)\|_{L^2(\Omega)}^2+\|\nabla
n_k(t)\|_{L^2(\Omega)}^2\big)\,dt\le
C(\beta,\delta,\rho_0,n_0,M_0),\\[4mm]
\|\rho_k+n_k\|_{L^{\beta+1}(Q_T)}\le
C(\epsilon,\beta,\delta,\rho_0,n_0,M_0), \end{cases}\ee where
$Q_T=\Omega\times (0,T)$ and $\beta\ge4$.

This yields the following Proposition by the analysis in
\cite{Feireisl} (see also \cite{Vasseur-Wen-Yu}).
\begin{proposition}\label{0-le:aweak solution}
Suppose $\beta>\max\{4,\Gamma+1,\gamma+1\}$. For any given
$\epsilon,\delta>0$, there exists a global weak solution
$(\rho,n,u)$ to (\ref{a-equation}), (\ref{a-initial}) and
(\ref{a-boundary}) such that for any given $T>0$, the following
estimates
 \be\label{1-r1}
\sup\limits_{t\in[0,T]}\|(\rho_+,\rho)(t)\|_{L^\gamma(\Omega)}^\gamma\le
C(\rho_0,n_0,M_0),\ee \be\label{1-r1+1}
\sup\limits_{t\in[0,T]}\|(\rho_-,n)(t)\|_{L^\Gamma(\Omega)}^\Gamma\le
C(\rho_0,n_0,M_0),\ee
\be\label{1-r2}\delta\sup\limits_{t\in[0,T]}\|(\rho,n)(t)\|_{L^\beta(\Omega)}^\beta\le
C(\rho_0,n_0,M_0),\ee
\be\label{1-r3}\sup\limits_{t\in[0,T]}\|\sqrt{\rho+n}u(t)\|_{L^2(\Omega)}^2\le
C(\rho_0,n_0,M_0),\ee \be\label{1-r4}
\int_0^T\|u(t)\|_{H^1_0(\Omega)}^2\,dt\le C(\rho_0,n_0,M_0),\ee
\be\label{1-r5} \epsilon\int_0^T\|(\nabla\rho,\nabla
n)(t)\|_{L^2(\Omega)}^2\,dt\le C(\beta,\delta,\rho_0,n_0,M_0),\ee
and \be\label{1-r6} \|(\rho,n)(t)\|_{L^{\beta+1}(Q_T)}\le
C(\epsilon,\beta,\delta,\rho_0,n_0,M_0) \ee hold, where the norm
$\|(\cdot,\cdot)\|$ denotes $\|\cdot\| + \|\cdot\|$, and
$\rho,n\ge0$ a.e. on $Q_T$.

 Finally, there exists $r>1$ such
that $\rho_t, n_t, \nabla^2\rho, \nabla^2 n\in L^r(Q_T)$ and the
equations (\ref{a-equation})$_1$ and (\ref{a-equation})$_2$ are
satisfied a.e. on $Q_T$.
\end{proposition}


\section{The vanishing of the artificial viscosity}\label{sec4}
\setcounter{equation}{0} \setcounter{theorem}{0}

 In this section, let $C$ denote a generic positive constant
depending on the initial data, $\delta$ and some other known
constants but independent of $\epsilon$.

\subsection{Passing to the limit as
$\epsilon\rightarrow0^+$}\label{s4.1}

The uniform estimates for $\epsilon$ resulting from (\ref{1-r1}),
(\ref{1-r1+1}), and (\ref{1-r2}) are not sufficient to obtain the weak
convergence of the artificial pressure
$P(n_{\epsilon},\rho_{\epsilon})+\delta(\rho_\epsilon+n_\epsilon)^\beta$
which is bounded only in $L^1(Q_T)$. Thus we need to obtain higher
integrability estimate of the artificial pressure uniformly for
$\epsilon$.

In the rest of the section, we remove the subscript $\epsilon$ of
the solutions for brevity.
\begin{lemma}\label{le:h-inofrho}
Let $(\rho,n,u)$ be the solution given by Proposition
\ref{0-le:aweak solution}, then \bex \int_0^T\int_\Omega
(n^{\Gamma+1}+\rho^{\gamma+1}+\delta\rho^{\beta+1}+\delta
n^{\beta+1})\,dx\,dt\le C\eex for
$\beta>\max\{4,\Gamma+1,\gamma+1\}$.
\end{lemma}
\pf The proof can be done by using (\ref{pressure ineq}) and the
arguments similar to \cite{Feireisl} where the test function
$\psi(t)\mathcal{B}[\rho-\widehat{\rho}]$ is replaced by
$\psi(t)\mathcal{B}[\rho+n-\widehat{\rho+n}]$. Here
\bex\begin{split}&\mathcal{B}: \Big\{f\in L^p(\Omega);\,
|\Omega|^{-1}\int_\Omega f\,dx=0\Big\}\mapsto
W_{0}^{1,p}(\Omega),\quad 1<p<\infty,\\&\psi\in
C_0^\infty(0,T),\,\, 0\le \psi\le 1,\, \mathrm{and}\,\,
\widehat{G}=\frac{1}{|\Omega|}\int_\Omega G \,dx\end{split}\eex
for $G=\rho,\,\rho+n$.
\endpf

Due to the relation between $P$ and $(n,\rho)$, i.e.,
(\ref{pressure ineq}), we have the following corollary.

\begin{corollary}\label{cor4.2}
Let $(\rho,n,u)$ be the solution given by Proposition
\ref{0-le:aweak solution}, then \bex \int_0^T\int_\Omega
\big({\rho_+}^{\gamma_1}+{\rho_-}^{\Gamma_1}\big)\,dx\,dt\le C,
\eex where
$\gamma_1=\gamma\min\big\{\frac{\gamma+1}{\gamma},\frac{\Gamma+1}{\Gamma}\big\}$
and
$\Gamma_1=\Gamma\min\big\{\frac{\gamma+1}{\gamma},\frac{\Gamma+1}{\Gamma}\big\}$.
Note that $\gamma_1>\gamma$ and $\Gamma_1>\Gamma$.

\end{corollary}

\bigskip

With (\ref{1-r1})-(\ref{1-r5}) and Lemma \ref{le:h-inofrho} and
Corollary \ref{cor4.2}, we are able to pass to the limits as
$\epsilon\rightarrow0^+$. Before doing this, we need to dress the
approximate solution constructed in Proposition \ref{0-le:aweak
solution} in the lower subscript ``$\epsilon$" for fixed
$\delta>0$, i.e., $(\rho_\epsilon,n_\epsilon, u_\epsilon)$. Then
letting $\epsilon\rightarrow0^+$ (taking a subsequence if
necessary), we have \be\label{2-lim}
\begin{cases}
(\rho_\epsilon,n_\epsilon)\rightarrow (\rho,n)\ \mathrm{in}\
C([0,T]; L_{weak}^\beta(\Omega))\cap C([0,T];H^{-1}(\Omega))\
\mathrm{and}\ \mathrm{weakly}\
\mathrm{in}\ L^{\beta+1}(Q_T),\\[2mm]
(\epsilon\Delta\rho_\epsilon, \epsilon\Delta
n_\epsilon)\rightarrow 0\ \mathrm{weakly}\ \mathrm{in}\ L^2(0,T;
H^{-1}(\Omega)),\\[2mm]
u_\epsilon\rightarrow u\ \mathrm{weakly}\ \mathrm{in}\ L^2(0,T;
H_0^1(\Omega)),\\[2mm]
(\rho_\epsilon+n_\epsilon) u_\epsilon\rightarrow (\rho+n) u \
\mathrm{in}\ C([0,T]; L_{weak}^\frac{2\beta}{\beta+1})\cap
C([0,T];H^{-1}(\Omega)),\\[2mm]
(\rho_\epsilon u_\epsilon,n_\epsilon u_\epsilon)\rightarrow (\rho
u, nu) \ \mathrm{in}\ \mathcal{D}^\prime(Q_T),\\[2mm]
(\rho_\epsilon+n_\epsilon) u_\epsilon\otimes u_\epsilon\rightarrow
(\rho+n) u\otimes u\ \mathrm{in}\ \mathcal{D}^\prime(Q_T),\\[2mm]
P(n_\epsilon,\rho_\epsilon)+\delta(\rho_\epsilon+n_\epsilon)^\beta\rightarrow
\overline{P(n,\rho)+\delta(\rho+n)^\beta}\
\mathrm{weakly}\ \mathrm{in}\ L^\frac{\beta+1}{\beta}(Q_T),\\[2mm]
\epsilon\nabla
u_\epsilon\cdot\nabla(\rho_\epsilon+n_\epsilon)\rightarrow 0\
\mathrm{in}\ L^1(Q_T),
\end{cases}
\ee and $\rho,n\ge0$, where the limit $(\rho, n, u)$ solves the
following system in the sense of distribution on $Q_T$ for any
$T>0$:
 \begin{equation}\label{a2-equation}
    \left\{
    \begin{array}{l}
        n_t+\mathrm{div}(n u)=0,\\[2mm]
        \rho_t+\mathrm{div}(\rho u)=0,\\[2mm]
 \big[(\rho+n) u\big]_t+\mathrm{div}\big[(\rho+n) u\otimes u\big]
 +\nabla\overline{P(n,\rho)+\delta(\rho+n)^\beta} =
            \mu \Delta u + (\mu+\lambda)\nabla \mathrm{div}u
    \end{array}
    \right.
    \end{equation}
    with initial and boundary condition
\be\label{a2-initial} \big(\rho, n,(\rho+n)
u\big)|_{t=0}=(\rho_{0,\delta},n_{0,\delta},M_{0,\delta}), \ee
\be\label{a2-boundary} u|_{\partial\Omega}=0, \ee where
$\overline{f(t,x)}$ denotes the weak limit of $f_{\epsilon}(t,x)$
as $\epsilon\to0.$

To this end, we have to show that
$$\overline{P(n,\rho)+\delta(\rho+n)^\beta}=P(n,\rho)+\delta(\rho+n)^\beta.$$

\subsection{The weak limit of the pressure}\label{s4.2}
This part is similar to \cite{Vasseur-Wen-Yu}, where it
focuses on the more complicated pressure $P$, since the artificial
pressure term $\delta(\rho_{\epsilon}+n_{\epsilon})^\beta$ controls the possible oscillation for
$(\rho_\epsilon+n_\epsilon)^\Gamma$ and
$(\rho_\epsilon+n_\epsilon)^\gamma$ arising in one of the
decomposition terms of the pressure, i.e.,
$P(Ad_{\epsilon},Bd_{\epsilon})$ where
$d_\epsilon=\rho_\epsilon+n_\epsilon$,
$(A,B)=(\frac{n}{\rho+n},\frac{\rho}{\rho+n})$ if $\rho+n\neq0$,
and $0\le A, B\le 1$, $\Big(A(\rho+n),B(\rho+n)\Big)=(n,\rho)$.

\bigskip

{\noindent\bf Claim:}
\be\label{2-claim1}\overline{P(n,\rho)+\delta(\rho+n)^\beta}=P(n,\rho)+\delta(\rho+n)^\beta
\ee a.e. on $Q_T$.

\medskip

The proof of (\ref{2-claim1}) relies on the following lemmas. In
particular, the next lemma plays an essential role.

\begin{lemma}\label{le:2-pressure}Let $(\rho_\epsilon,n_\epsilon)$ be the solution stated in Proposition
\ref{0-le:aweak solution}, and $(\rho, n)$ be the limit in the
sense of (\ref{2-lim}), then \be\label{2-a1} (\rho+n)
\overline{P(n,\rho)} \le \overline{(\rho+n) P(n,\rho)}\ee a.e. on
$\Omega\times(0,T)$.
\end{lemma}

\pf The idea is similar to \cite{Vasseur-Wen-Yu} by Vasseur, the
author, and Yu. However, since the pressure here is more
complicated, we have to give a complete proof.

As in \cite{Vasseur-Wen-Yu}, the pressure and
$n_\epsilon+\rho_\epsilon$ are decomposed as
follows.\be\label{2-decom of press}
\begin{cases}
P(n_\epsilon,\rho_\epsilon)=P(A_\epsilon d_\epsilon, B_\epsilon
d_\epsilon)-P(A
d_\epsilon,B d_\epsilon) + P(A d_\epsilon,B d_\epsilon),\\[2mm]
n_\epsilon+\rho_\epsilon=(A_\epsilon+B_\epsilon) d_\epsilon
=(A+B)d_\epsilon +(A_\epsilon-A+B_\epsilon-B) d_\epsilon,
\end{cases}
\ee where  $d_\epsilon=\rho_\epsilon+n_\epsilon$, $d=\rho+n$,
$(A_\epsilon,B_\epsilon)=(\frac{n_\epsilon}{d_\epsilon},\frac{\rho_\epsilon}{d_\epsilon})$
if $d_\epsilon\neq0$, $(A,B)=(\frac{n}{d},\frac{\rho}{d})$ if
$d\neq0$, $0\le A_\epsilon,B_\epsilon, A, B\le 1$, and
$(A_\epsilon d_\epsilon,B_\epsilon
d_\epsilon)=(n_\epsilon,\rho_\epsilon), (Ad,Bd)=(n,\rho)$,
$(\rho,n)$ is the limit of $(\rho_{\epsilon},n_{\epsilon})$ in a
suitable weak topology.

For any $\psi\in C([0,t])$, $\phi\in C(\overline{\Omega})$ where
$\psi,\phi\ge0$, we use (\ref{2-decom of press}) and obtain
\be\label{2-II}\begin{split} &\int_0^t\psi\int_\Omega\phi
P(n_\epsilon,\rho_\epsilon)(\rho_\epsilon+n_\epsilon) \,dx\,ds\\=&
\int_0^t\psi\int_\Omega\phi P(A d_\epsilon, B d_\epsilon)
(A+B)d_\epsilon\,dx\,ds+\int_0^t\psi\int_\Omega\phi P(A
d_\epsilon, B d_\epsilon) (A_\epsilon-A+B_\epsilon-B)
d_\epsilon\,dx\,ds\\& +\int_0^t\psi\int_\Omega\phi
\Big[P(A_\epsilon d_\epsilon, B_\epsilon d_\epsilon)-P(A
d_\epsilon,B
d_\epsilon)\Big](\rho_\epsilon+n_\epsilon)\,dx\,dt\\=&
\sum\limits_{i=1}^3II_i.
\end{split}
\ee

 For $II_2$, we follow an argument similar to \cite{Vasseur-Wen-Yu}. More specifically, there exists a positive integer $k_0$ large
enough such that \be\label{2-18}
\max\{\frac{k_0\gamma}{k_0-1},\frac{k_0\Gamma}{k_0-1}\}\le\beta
\ee due to the assumption that $\max\{\Gamma,\gamma\}<\beta$.
Therefore (\ref{2-18}) implies that \be\label{3-1} \int_{Q_T}
d_\epsilon|d_\epsilon^\frac{k_0\Gamma}{k_0-1}+
d_\epsilon^\frac{k_0\gamma}{k_0-1}|\,dx\,dt\le C\int_{Q_T}
|d_\epsilon^{\beta+1}+1|\,dx\,dt\le C \ee where we have used Lemma \ref{le:h-inofrho} that $d_\epsilon$ is bounded in $L^{\beta+1}(Q_T)$
uniformly for $\epsilon$.

Recalling (\ref{pressure ineq}), we have \be\label{pressure
ineq+1} P(Ad_\epsilon,B d_\epsilon)\le C_0(A^\Gamma
d_\epsilon^{\Gamma}+B^\gamma d_\epsilon^{\gamma}). \ee This
together with H\"older inequality and (\ref{3-1}) yields
\be\label{2-II2}\begin{split} |II_2| \le&C\left(\int_{Q_T}
d_\epsilon|A_\epsilon-A|^{k_0}\,dx\,dt\right)^\frac{1}{k_0}\left(\int_{Q_T}
d_\epsilon|A^\Gamma d_\epsilon^\Gamma+B^\gamma
d_\epsilon^\gamma|^\frac{k_0}{k_0-1}\,dx\,dt\right)^\frac{k_0-1}{k_0}\\&+
C\left(\int_{Q_T}
d_\epsilon|B_\epsilon-B|^{k_0}\,dx\,dt\right)^\frac{1}{k_0}\left(\int_{Q_T}
d_\epsilon|A^\Gamma d_\epsilon^\Gamma+B^\gamma
d_\epsilon^\gamma|^\frac{k_0}{k_0-1}\,dx\,dt\right)^\frac{k_0-1}{k_0}
\\
\le&C\left(\int_{Q_T}
d_\epsilon|A_\epsilon-A|^{k_0}\,dx\,dt\right)^\frac{1}{k_0}\left(\int_{Q_T}
d_\epsilon|d_\epsilon^\frac{k_0\Gamma}{k_0-1}+
d_\epsilon^\frac{k_0\gamma}{k_0-1}|\,dx\,dt\right)^\frac{k_0-1}{k_0}\\&+
C\left(\int_{Q_T}
d_\epsilon|B_\epsilon-B|^{k_0}\,dx\,dt\right)^\frac{1}{k_0}\left(\int_{Q_T}
d_\epsilon|d_\epsilon^\frac{k_0\Gamma}{k_0-1}+
d_\epsilon^\frac{k_0\gamma}{k_0-1}|\,dx\,dt\right)^\frac{k_0-1}{k_0}
\\
\le&C\left(\int_{Q_T}
d_\epsilon|A_\epsilon-A|^{k_0}\,dx\,dt\right)^\frac{1}{k_0}+
C\left(\int_{Q_T}
d_\epsilon|B_\epsilon-B|^{k_0}\,dx\,dt\right)^\frac{1}{k_0}.
\end{split}
\ee Choosing $\nu_k:=\nu_\epsilon=\epsilon$ in Lemma \ref{main
2-le:important}, we conclude that
 \be\label{4.14}\begin{split}
 &\left(\int_{Q_T} d_\epsilon|A_\epsilon-A|^{k_0}\,dx\,dt\right)^\frac{1}{k_0}\to
 0,\\
 &\left(\int_{Q_T} d_\epsilon|B_\epsilon-B|^{k_0}\,dx\,dt\right)^\frac{1}{k_0}\to 0
 \end{split}
 \ee
as $\epsilon$ goes to zero. In fact, $d_{\epsilon}\in
L^{\infty}(0,T;L^{\beta}(\Omega))$
 for $\beta>4$, and $u_{\epsilon}\in L^2(0,T;H^1_0(\Omega)),$ and
 \begin{equation*}
  \sqrt{\epsilon}\|\nabla\rho_{\epsilon}\|_{L^2(0,T;L^2(\Omega))}\leq C_0,\,\sqrt{\epsilon}\|\nabla n_{\epsilon}\|_{L^2(0,T;L^2(\Omega))}
  \leq C_0,
\end{equation*}
 and for any $\epsilon>0$ and any $t> 0$:
\begin{equation}
\label{initial condition for n2 over density}\int_{\Omega}
\frac{b_{\epsilon}^2}{d_{\epsilon}}\,dx\leq
\int_{\Omega}\frac{b_0^2}{d_0}\,dx
\end{equation} where $d_\epsilon=\rho_\epsilon+n_\epsilon$, $b_\epsilon=\rho_\epsilon$,
$n_\epsilon$, and (\ref{initial condition for n2 over density}) is
obtained  in Remark 2.4, \cite{Vasseur-Wen-Yu}. Thus, we are able
to apply Lemma \ref{main 2-le:important} to deduce \eqref{4.14}.
Hence we have $II_2\to 0$ as $\epsilon\to 0$.

For $II_3$, the analysis becomes more complicated due to the
pressure. First, we need the following estimate.
\be\label{pressurediff}
\begin{split}
P(A_\epsilon d_\epsilon, B_\epsilon d_\epsilon)-P(A d_\epsilon,B
d_\epsilon)=&A_+\gamma[\rho_+(\xi_1,\xi_2)]^{\gamma-1}\partial_{\xi_1}\rho_+(\xi_1,\xi_2)\big[A_\epsilon
d_\epsilon-Ad_\epsilon\big]\\&+A_+\gamma[\rho_+(\xi_1,\xi_2)]^{\gamma-1}\partial_{\xi_2}\rho_+(\xi_1,\xi_2)\big[B_\epsilon
d_\epsilon-Bd_\epsilon\big]\\=&\frac{A_+\gamma\big(\frac{A_-}{A_+}\big)^\frac{1}{\Gamma}
[\rho_+(\xi_1,\xi_2)]^{\gamma-\frac{\gamma}{\Gamma}}}{\frac{\gamma}{\Gamma}
[1-\alpha(\xi_1,\xi_2)]
 +\alpha(\xi_1,\xi_2)}\big[A_\epsilon
d_\epsilon-Ad_\epsilon\big]\\&+\frac{A_+\gamma[\rho_+(\xi_1,\xi_2)]^{\gamma-1}}{\frac{\gamma}{\Gamma}[1-\alpha(\xi_1,\xi_2)]
 +\alpha(\xi_1,\xi_2)}\big[B_\epsilon
d_\epsilon-Bd_\epsilon\big],
\end{split}
\ee since \be\label{partialrhon} \begin{cases}
\frac{\partial\rho_+(n,\rho)}{\partial
n}=\frac{\big(\frac{A_-}{A_+}\big)^\frac{1}{\Gamma}\rho_+^{1-\frac{\gamma}{\Gamma}}}{\frac{\gamma}{\Gamma}(1-\alpha)
 +\alpha},\\[2mm]
\frac{\partial\rho_+(n,\rho)}{\partial\rho}=\frac{1}{\frac{\gamma}{\Gamma}(1-\alpha)
 +\alpha}
\end{cases}
 \ee which can be obtained similarly to (\ref{narho+}). Here $\xi_1$ $(\xi_2)$
 varies between $A_\epsilon d_\epsilon$ $(B_\epsilon d_\epsilon)$
 and $A d_\epsilon$ ($B d_\epsilon$).

 In view of (\ref{pressure ineq}), we have
 \be\label{rho+bound}
 \begin{split}
\rho_+(\xi_1,\xi_2)\le
C\Big(\xi_1^\frac{\Gamma}{\gamma}+\xi_2\Big)\le
C\Big[(\rho_\epsilon+n_\epsilon)^\frac{\Gamma}{\gamma}+\rho_\epsilon+n_\epsilon\Big],
 \end{split}
 \ee where we have used
 \bex
0\le\xi_1,\xi_2\le \rho_\epsilon+n_\epsilon.
 \eex

 By virtue of (\ref{pressurediff}) and (\ref{rho+bound}), and using Young inequality, we have
 \be\label{Pressurezero}\begin{split}
&|P(A_\epsilon d_\epsilon, B_\epsilon d_\epsilon)-P(A d_\epsilon,B
d_\epsilon)|\\
\le&C[(\rho_\epsilon+n_\epsilon)^\frac{\Gamma}{\gamma}+\rho_\epsilon+n_\epsilon]^{\gamma-\frac{\gamma}{\Gamma}}\big|A_\epsilon
d_\epsilon-Ad_\epsilon\big|+C[(\rho_\epsilon+n_\epsilon)^\frac{\Gamma}{\gamma}+\rho_\epsilon+n_\epsilon]^{\gamma-1}\big|B_\epsilon
d_\epsilon-Bd_\epsilon\big|
\\
\le&C\Big[(\rho_\epsilon+n_\epsilon)^\Gamma+(\rho_\epsilon+n_\epsilon)^\gamma+1\Big]\Big(\big|A_\epsilon
d_\epsilon-Ad_\epsilon\big|+\big|B_\epsilon
d_\epsilon-Bd_\epsilon\big|\Big)\\
\le&C\Big[d_\epsilon^{\Gamma_m}+1\Big]\Big(\big|A_\epsilon
d_\epsilon-Ad_\epsilon\big|+\big|B_\epsilon
d_\epsilon-Bd_\epsilon\big|\Big),
 \end{split}
 \ee where $\Gamma_m=\max\big\{\Gamma,\gamma\big\}$.

 Now we are in a position to estimate $II_3$. In fact, there exists a positive integer $k_1$ large
enough such that \be\label{2-au}\begin{split}
(\Gamma_m+2-\frac{1}{k_1})\frac{k_1}{k_1-1}<\beta+1
\end{split}\ee due to the assumption $\gamma+1,\Gamma+1<\beta$.

In virtue of (\ref{Pressurezero}), we have\bex
\begin{split}
 |II_3|=&\int_0^t\psi\int_\Omega\phi
\Big[P(A_\epsilon d_\epsilon, B_\epsilon d_\epsilon)-P(A
d_\epsilon,B d_\epsilon)\Big](\rho_\epsilon+n_\epsilon)\,dx\,dt\\
\le& C\int_{Q_T}\big(d_\epsilon^{\Gamma_m+1}+1\big)\big|A_\epsilon
d_\epsilon-Ad_\epsilon\big|\,dx\,ds+C\int_{Q_T}\big(d_\epsilon^{\Gamma_m+1}+1\big)\big|B_\epsilon
d_\epsilon-Bd_\epsilon\big|\,dx\,ds\\
=&
C\int_{Q_T}\big(d_\epsilon^{\Gamma_m+2-\frac{1}{k_1}}d_\epsilon^\frac{1}{k_1}+d_\epsilon^\frac{1}{2}d_\epsilon^\frac{1}{2}\big)\big|A_\epsilon-A\big|\,dx\,ds
+C\int_{Q_T}\big(d_\epsilon^{\Gamma_m+2-\frac{1}{k_1}}d_\epsilon^\frac{1}{k_1}+d_\epsilon^\frac{1}{2}d_\epsilon^\frac{1}{2}\big)\big|B_\epsilon
-B\big|\,dx\,ds.
\end{split}
\eex Then applying H\"older inequality, we get
\be\label{2-II3}\begin{split}
 |II_3|\le&
C\left(\int_{Q_T}
d_\epsilon^{(\Gamma_m+2-\frac{1}{k_1})\frac{k_1}{k_1-1}}\,dx\,dt\right)^\frac{k_1-1}{k_1}
\left(\int_{Q_T} d_\epsilon
\big|A_\epsilon-A\big|^{k_1}\,dx\,dt\right)^\frac{1}{k_1}\\&+C\left(\int_{Q_T}
d_\epsilon\,dx\,dt\right)^\frac{1}{2} \left(\int_{Q_T} d_\epsilon
\big|A_\epsilon-A\big|^{2}\,dx\,dt\right)^\frac{1}{2}\\&+
C\left(\int_{Q_T}
d_\epsilon^{(\Gamma_m+2-\frac{1}{k_1})\frac{k_1}{k_1-1}}\,dx\,dt\right)^\frac{k_1-1}{k_1}
\left(\int_{Q_T} d_\epsilon
\big|B_\epsilon-B\big|^{k_1}\,dx\,dt\right)^\frac{1}{k_1}\\&+C\left(\int_{Q_T}
d_\epsilon\,dx\,dt\right)^\frac{1}{2} \left(\int_{Q_T} d_\epsilon
\big|B_\epsilon-B\big|^{2}\,dx\,dt\right)^\frac{1}{2}
\\ \le& C\left(\int_{Q_T} d_\epsilon
\big|A_\epsilon-A\big|^{k_1}\,dx\,dt\right)^\frac{1}{k_1}+C\left(\int_{Q_T}
d_\epsilon
\big|A_\epsilon-A\big|^{2}\,dx\,dt\right)^\frac{1}{2}\\&+
C\left(\int_{Q_T} d_\epsilon
\big|B_\epsilon-B\big|^{k_1}\,dx\,dt\right)^\frac{1}{k_1} +
C\left(\int_{Q_T} d_\epsilon
\big|B_\epsilon-B\big|^{2}\,dx\,dt\right)^\frac{1}{2}\rightarrow 0
\end{split}
\ee as $\epsilon\rightarrow0^+$, where we have used
 (\ref{4.14}), (\ref{2-au}), Lemma \ref{le:h-inofrho}, and Young inequality.

Combining (\ref{2-II}), (\ref{2-II2}) and (\ref{2-II3}), we have
\be\label{2-II+1}\begin{split}
\lim\limits_{\epsilon\rightarrow0^+}\int_0^t\psi\int_\Omega\phi
P(n_\epsilon,\rho_\epsilon)(\rho_\epsilon+n_\epsilon) \,dx\,ds=&
\int_0^t\psi\int_\Omega\phi(A+B)\overline{\overline{P(A d, B
d)d}}\,dx\,ds\\
\ge& \int_0^t\psi\int_\Omega\phi(A+B)\overline{\overline{P(A d, B
d)}}d\,dx\,ds\\
=& \int_0^t\psi\int_\Omega\phi(\rho+n)\overline{\overline{P(A d, B
d)}}\,dx\,ds
\end{split}
\ee  where we have used that $A+B=1$, and Lemma
\ref{Feireisl-Novotny lemma} such that \bex\overline{\overline{P(A
d, B d)d}}\ge \overline{\overline{P(A d, B d)}}d\eex due to the
fact that the functions $z\mapsto P(A z, B z)$ and $z\mapsto z$
are non decreasing functions. Here $\overline{\overline{(\cdot)}}$
represents the weak limit of $(\cdot)$ with respect to
$d_\epsilon$ as $\epsilon\rightarrow0^+$. Note that in this
section $P(Ad_\epsilon,Bd_\epsilon)d_\epsilon$ and
$P(Ad_\epsilon,Bd_\epsilon)$ are bounded in
$L^\frac{\beta+1}{\Gamma_m+1}(Q_T)$ and in
$L^\frac{\beta+1}{\Gamma_m}(Q_T)$, respectively, due to Lemma
\ref{le:h-inofrho}. Moreover, both $\frac{\beta+1}{\Gamma_m+1}$
and $\frac{\beta+1}{\Gamma_m}$
 are large than 1, which implies that
$\overline{\overline{P(A d, B d)d}}$ and $\overline{\overline{P(A
d, B d)}}$ are well-defined.

We claim that \be\label{2-claim}
\int_0^t\psi\int_\Omega\phi(\rho+n)\overline{\overline{P(A d, B
d)}}\,dx\,ds=\int_0^t\psi\int_\Omega\phi(\rho+n)\overline{P(n,
\rho)}\,dx\,ds. \ee In fact,  \be\label{2-419}
\begin{split}
&\int_0^t\psi\int_\Omega\phi (\rho+n)\overline{\overline{P(A d, B
d)}}\,dx\,ds
\\=&\lim\limits_{\epsilon\rightarrow0^+}\int_0^t\psi\int_\Omega\phi
(\rho+n) P(A d_\epsilon, B d_\epsilon)\,dx\,ds
\\=&\lim\limits_{\epsilon\rightarrow0^+}\int_0^t\psi\int_\Omega\phi
(\rho+n) P(n_\epsilon,
\rho_\epsilon)\,dx\,ds+\lim\limits_{\epsilon\rightarrow0^+}\int_0^t\psi\int_\Omega\phi
(\rho+n) \Big[P(A d_\epsilon, B d_\epsilon)-P(n_\epsilon,
\rho_\epsilon)\Big]\,dx\,ds\\=&\int_0^t\psi\int_\Omega\phi
(\rho+n) \overline{P(n,
\rho)}\,dx\,ds+\lim\limits_{\epsilon\rightarrow0^+}\int_0^t\psi\int_\Omega\phi
(\rho+n) \Big[P(A d_\epsilon, B d_\epsilon)-P(A_\epsilon
d_\epsilon, B_\epsilon d_\epsilon)\Big]\,dx\,ds.
\end{split}
\ee Similar to $II_3$, the last term on the right hand side of
(\ref{2-419}) converges to zero as $\epsilon\rightarrow0^+$. Hence
we get (\ref{2-claim}).

In view of (\ref{2-II+1}), (\ref{2-claim}), and the fact that the
test functions $\phi$ and $\psi$ are arbitrary, we complete the
proof of the lemma.

\endpf

\begin{lemma}\label{2-le:3.7} Let
$(\rho_\epsilon,n_\epsilon,u_\epsilon)$ be the solution stated in
Proposition \ref{0-le:aweak solution}, and $(\rho, n,u)$ be the
limit in the sense of (\ref{2-lim}), then \be\label{2-8}
\lim\limits_{\epsilon\rightarrow0^+}\int_{Q_T}\psi\phi
H_\epsilon(\rho_\epsilon+n_\epsilon)\,dx\,dt=\int_{Q_T}\psi\phi
H(\rho+n)\,dx\,dt, \ee for any $\psi\in  C_0^\infty(0,T)$ and
$\phi\in
 C_0^\infty(\Omega)$, where  \bex\begin{split}
H_\epsilon:=&P(n_\epsilon,\rho_\epsilon)+\delta(\rho_\epsilon+n_\epsilon)^\beta-(2\mu+\lambda)\mathrm{div}u_\epsilon,
\\
H:=&\overline{P(n,\rho)+\delta(\rho+n)^\beta}-(2\mu+\lambda)\mathrm{div}u.
\end{split}
\eex
\end{lemma}
\begin{remark}
 The proof of (\ref{2-8}) is motivated by \cite{Feireisl} for Navier-Stokes equations.
In fact, the lemma can be found in \cite{Vasseur-Wen-Yu} where the
pressure is given by (\ref{pressure law1}). For the pressure
(\ref{pressure law2}), the proof is similar.
\end{remark}

\medskip

With Lemmas \ref{le:2-pressure} and \ref{2-le:3.7}, it is not
difficult to obtain the next lemma.

\begin{lemma}Let $(\rho_\epsilon,n_\epsilon)$ be the solution stated in Lemma
\ref{0-le:aweak solution}, and $(\rho, n)$ be the limit in the
sense of (\ref{2-lim}), then \be\label{2-a2}\int_0^t\int_\Omega
(\rho+n)
\mathrm{div}u\,dx\,ds\le\lim\limits_{\epsilon\rightarrow0^+}\int_0^t\int_\Omega
(\rho_\epsilon+n_\epsilon) \mathrm{div}u_\epsilon\,dx\,ds \ee for
a.e. $t\in(0,T)$.
\end{lemma}

\medskip

By virtue of Lemma 4.4 in \cite{Vasseur-Wen-Yu}, we have
 \be\label{2-17}
\begin{split}
&\int_\Omega \big[\rho_\epsilon\log \rho_\epsilon-\rho\log
\rho+n_\epsilon\log n_\epsilon-n\log n\big](t)\,dx
 \\ \le&\int_0^t\int_\Omega
(\rho+n) \mathrm{div}u\,dx\,ds-\int_0^t\int_\Omega
(\rho_\epsilon+n_\epsilon) \mathrm{div}u_\epsilon\,dx\,ds
\end{split}
\ee for a.e. $t\in(0,T)$.

Passing both sides of (\ref{2-17}) to the limits as
$\epsilon\rightarrow0^+$, and using (\ref{2-a2}), we have \bex
\begin{split}
\int_\Omega \big[\overline{\rho\log \rho}-\rho\log
\rho+\overline{n\log n}-n\log n\big](t)\,dx \le0.
\end{split}
\eex Thanks to the convexity of $z\mapsto z\log z$, we have \bex
\overline{\rho\log \rho}\ge\rho\log \rho\quad \mathrm{and}\quad
\overline{n\log n}\ge n \log n \eex a.e. on $Q_T$. This turns out
that \bex \int_\Omega \big[\overline{\rho\log \rho}-\rho\log
\rho+\overline{n\log n}-n\log n\big](t)\,dx=0. \eex Hence we get
\bex \overline{\rho\log \rho}=\rho\log \rho\quad \mathrm{and}\quad
\overline{n\log n}= n \log n \eex a.e. on $Q_T$, which implies
that $(\rho_\epsilon,n_\epsilon)\rightarrow (\rho,n)$ a.e. in
$Q_T$. It combined with Lemma \ref{le:h-inofrho} yields strong
convergence of $(\rho_\epsilon,n_\epsilon)$ in $L^{\beta_1}(Q_T)$
for any $\beta_1<\beta+1$. Thus we complete the proof of
(\ref{2-claim1}).

\bigskip

To this end, we give a proposition as a summary for this section.

\begin{proposition}\label{2-le:aweak solution}
Suppose $\beta>\max\{4,\Gamma+1,\gamma+1\}$. For any given
$\delta>0$, there exists a global weak solution
$(\rho_\delta,n_\delta,u_\delta)$ to the following system over
$\Omega\times(0,\infty)$:
 \begin{equation}\label{a3-equation}
    \left\{
    \begin{array}{l}
        n_t+\mathrm{div}(n u)=0,\\[2mm]
        \rho_t+\mathrm{div}(\rho u)=0,\\[2mm]
 \big[(\rho+n) u\big]_t+\mathrm{div}\big[(\rho+n) u\otimes u\big]
 +\nabla P(n,\rho)+\delta\nabla (\rho+n)^\beta =
            \mu \Delta u + (\mu+\lambda)\nabla \mathrm{div}u,
    \end{array}
    \right.
    \end{equation}
    with initial and boundary condition
\be\label{a3-initial}\big (\rho, n,(\rho+n)
u\big)|_{t=0}=(\rho_{0,\delta},n_{0,\delta},M_{0,\delta})\
\mathrm{on}\ \overline{\Omega}, \ee \be\label{a3-boundary}
u|_{\partial\Omega}=0 \ \mathrm{for}\ t\ge0, \ee such that for any
given $T>0$, the following estimates
 \be\label{2-r1}
\sup\limits_{t\in[0,T]}\|\rho_\delta(t)\|_{L^\gamma(\Omega)}^\gamma\le
C(\rho_0,n_0,M_0),\ee

\be\label{2-r1+1}
\sup\limits_{t\in[0,T]}\|n_\delta(t)\|_{L^\Gamma(\Omega)}^\Gamma\le
C(\rho_0,n_0,M_0),\ee

\be\label{2-r2}\delta\sup\limits_{t\in[0,T]}\|(\rho_\delta(t),n_\delta(t))\|_{L^\beta(\Omega)}^\beta\le
C(\rho_0,n_0,M_0),\ee

\be\label{2-r3}\sup\limits_{t\in[0,T]}\|\sqrt{\rho_\delta+n_\delta}(t)u_\delta(t)\|_{L^2(\Omega)}^2\le
C(\rho_0,n_0,M_0),\ee

\be\label{2-r4} \int_0^T\|u_\delta(t)\|_{H^1_0(\Omega)}^2\,dt\le
C(\rho_0,n_0,M_0),\ee and \be\label{2-r5}
\|(\rho_\delta(t),n_\delta(t))\|_{L^{\beta+1}(Q_T)}\le
C(\beta,\delta,\rho_0,n_0,M_0) \ee hold, where the norm
$\|(\cdot,\cdot)\|$ denotes $\|\cdot\| + \|\cdot\|$.
\end{proposition}

\section{The vanishing of the artificial pressure}\label{sec5} \setcounter{equation}{0}
\setcounter{theorem}{0}

 Let $C$ be a generic constant depending only on the initial data and some other known constants but independent of $\delta$, which will
be used throughout this section.

\subsection{Passing to the limit as $\delta\rightarrow0^+$}

In this section,  we will obtain the global existence of the weak
solution to \eqref{equation}-\eqref{initial-boundary} by passing to the
limit of $(\rho_\delta,n_\delta,u_\delta)$ as
$\delta\rightarrow0^+$. To begin with, we have to get the higher
integrability estimates of the pressure $P$ uniformly for $\delta$
for the same reason as in the previous section.


In fact, as in \cite{Feireisl} (see
also \cite{Vasseur-Wen-Yu}), we have the following lemma.
\begin{lemma}\label{3-le:h-inofrho} Let $(\rho_\delta,n_\delta,u_\delta)$ be the solution
stated in Proposition \ref{2-le:aweak solution}, then we have
\be\label{3-hierho} \int_{Q_T}
(n_\delta^{\Gamma+\theta_1}+\rho_\delta^{\gamma+\theta_2}+\delta
n_\delta^{\beta+\theta_1}+\delta\rho_\delta^{\beta+\theta_2})\,dx
dt\le C(\theta_1,\theta_2)\ee for any positive constants
$\theta_1$ and $\theta_2$ satisfying
\bex\begin{split}\theta_1<\frac{\Gamma}{3}\,\, \mathrm{and}\,\,
\theta_1\le\min\{1,\frac{2\Gamma}{3}-1\};\,\,
\theta_2<\frac{\gamma}{3}\,\, \mathrm{and}\,\,
\theta_2\le\min\{1,\frac{2\gamma}{3}-1\}\quad \mathrm{if}\,\,
\Gamma,\gamma\in (\frac{3}{2},\infty).
\end{split}\eex
\end{lemma}

\bigskip

With (\ref{2-r1}), (\ref{2-r1+1}), (\ref{2-r3}), (\ref{2-r4}), and
(\ref{3-hierho}), letting $\delta\rightarrow0^+$ (taking a
subsequence if necessary), we have
 \be\label{3-lim}
\begin{cases}
\rho_\delta\rightarrow \rho\ \mathrm{in}\ C([0,T];
L_{weak}^\gamma(\Omega))\ \mathrm{and}\ \mathrm{weakly}\
\mathrm{in}\ L^{\gamma+\theta_2}(Q_T)\ \mathrm{as}\
\delta\rightarrow0^+,\\[2mm]
n_\delta\rightarrow n\ \mathrm{in}\ C([0,T];
L_{weak}^\Gamma(\Omega))\ \ \mathrm{and}\ \mathrm{weakly}\
\mathrm{in}\ L^{\Gamma+\theta_1}(Q_T)\ \mathrm{as}\
\delta\rightarrow0^+,\\[2mm]
u_\delta\rightarrow u\ \mathrm{weakly}\ \mathrm{in}\ L^2(0,T;
H_0^1(\Omega))\ \mathrm{as}\ \delta\rightarrow0^+,\\[2mm]
(\rho_\delta+n_\delta) u_\delta\rightarrow (\rho+n) u \
\mathrm{in}\ C([0,T];
L_{weak}^\frac{2\min\{\gamma,\Gamma\}}{\min\{\gamma,\Gamma\}+1})\cap
C([0,T];H^{-1}(\Omega))\ \mathrm{as}\ \delta\rightarrow0^+,\\[2mm]
(\rho_\delta u_\delta, n_\delta u_\delta)\rightarrow (\rho u, n u)
\ \mathrm{in}\ \mathcal{D}^\prime(Q_T)\ \mathrm{as}\
\delta\rightarrow0^+,\\[2mm]
(\rho_\delta+n_\delta) u_\delta\otimes u_\delta\rightarrow
(\rho+n) u\otimes u\ \mathrm{in}\ \mathcal{D}^\prime(Q_T)\
\mathrm{as}\
\delta\rightarrow0^+,\\[2mm]
P(n_\delta,\rho_\delta)\rightarrow \overline{P(n,\rho)}\
\mathrm{weakly}\ \mathrm{in}\
L^{\min\{\frac{\gamma+\theta_2}{\gamma},\frac{\Gamma+\theta_1}{\Gamma}\}}(Q_T)\
\mathrm{as}\
\delta\rightarrow0^+,\\[2mm]
\delta(\rho_\delta+n_\delta)^\beta\rightarrow 0\ \mathrm{in}\
L^1(Q_T)\ \mathrm{as}\ \delta\rightarrow0^+,
\end{cases}
\ee where the limit $(\rho, n, u)$ solves the following system in
the sense of distribution over $\Omega\times [0, T]$ for any given
$T>0$:
 \begin{equation}\label{3-equation}
    \left\{
    \begin{array}{l}
        n_t+\mathrm{div}(n u)=0,\\[2mm]
        \rho_t+\mathrm{div}(\rho u)=0,\\[2mm]
 \big[(\rho+n) u\big]_t+\mathrm{div}\big[(\rho+n) u\otimes u\big]
 +\nabla\overline{P(\rho,n)} =
            \mu \Delta u + (\mu+\lambda)\nabla \mathrm{div}u,
    \end{array}
    \right.
    \end{equation}
    with initial and boundary condition
\be\label{3-initial} (\rho, n,(\rho+n)
u)|_{t=0}=(\rho_0,n_0,M_0)\quad \mathrm{on}\ \overline{\Omega},
\ee \be\label{3-boundary} u|_{\partial\Omega}=0 \quad
\mathrm{for}\ t\ge0. \ee 

Finally, we need to justify that $\overline{P(\rho,n)}=P(\rho,n)$.
In fact, this has already been done by Vasseur, the author,
and Yu in \cite{Vasseur-Wen-Yu} for the pressure law
(\ref{pressure law1}) subject to the constraints
\be\label{additionalcondition}\max\{\frac{3\gamma}{4},\gamma-1,\frac{3(\gamma+1)}{5}\}<\Gamma<\min\{\frac{4\gamma}{3},\gamma+1,\frac{5\gamma}{3}-1\}\ee
and $\Gamma,\gamma>\frac{9}{5}$, which implies that $\Gamma$ and
$\gamma$ have to stay not too far from each other. Thus to
consider the case that $\Gamma,\gamma\ge\frac{9}{5}$ without any
other constraints, some new ingredients will be needed in the
following analysis.

\subsection{The weak limit of the pressure}\label{5.2}

To obtain the global existence of weak solution to
\eqref{equation}-\eqref{initial-boundary}, we have to justify the
following claim.

\bigskip

{\noindent\bf\large Claim.}
\be\label{3-claim}\overline{P(n,\rho)}=P(n,\rho)\ee for any
$\Gamma,\gamma\ge\frac{9}{5}$.

To prove (\ref{3-claim}), it suffices to derive the strong
convergence of $\rho_{\delta}$ and $n_{\delta}$ as
$\delta\rightarrow0^+$. In this section, we need that
$\rho_{\delta}$ and $n_{\delta}$ are bounded in $L^2(Q_T)$ for
that it will be essential to employ Lemma \ref{main
2-le:important}. As a consequence, the restriction that
$\gamma,\Gamma\ge\frac{9}{5}$ is needed in view of Lemma
\ref{3-le:h-inofrho}.

Lemmas \ref{3-hierho} and \ref{le:h-inofrho} indicate that the
uniform integrability of $\rho_\delta$ and $n_\delta$ is weaker
when $\Gamma,\gamma<3$. Thus some estimates such as (\ref{2-a1})
can not be obtained in this part. For this reason, we consider a
family of cut-off functions introduced in \cite{Feireisl} and
references therein, i.e., \be\label{Tk} T_k(z)=kT(\frac{z}{k}),\
z\in\mathbb{R},\ k=1,2,\cdot\cdot\cdot \ee where $T\in
C^\infty(\mathbb{R})$ satisfies \bex T(z)=\left\{\begin{array}{l} z\ for\ z\le1, \\
[3mm] 2\ for\ z\ge3,
\end{array}
\right.\eex and $T$ is concave.

The first conclusion in this subsection plays a very important
role, which is only subject to the constraint
$\Gamma,\gamma\ge\frac{9}{5}$. 
\begin{lemma}\label{3-le:4.6}
Let $(\rho_\delta,n_\delta)$ be the solutions constructed in
Proposition \ref{2-le:aweak solution}, and $(\rho,n)$ be the
limit, then  \be\label{3-c1}\begin{cases}
\overline{T_k(\rho)}\,\,\overline{P(n,\rho)}
 \le
\overline{T_k(\rho)P(n,\rho)},\\[4mm]
\overline{T_k(n)}\,\,\overline{P(n,\rho)}
 \le
\overline{T_k(n)P(n,\rho)},
\end{cases}
\ee a.e. on $\Omega\times(0,T)$, for any
$\Gamma,\gamma\ge\frac{9}{5}$.
\end{lemma}
\pf In view of Lemma \ref{main 2-le:important} with $\nu_K=0$ and $s=1$
$\Big($see \eqref{key bound}$\Big)$ where the condition (\ref{AAAinitial condition for n2 over
density}) can be ensured by using Lemma \ref{main lemma} for $N=2$ and $R=(n_\delta, \rho_\delta+n_\delta),$ $(\rho_\delta, \rho_\delta+n_\delta)$, we have
\be\label{ndconvergence}\begin{cases}
n_\delta-Ad_\delta\rightarrow0 \quad \mathrm{a.e.}\ \mathrm{in}\
Q_T,\\[2mm]\rho_\delta-Bd_\delta\rightarrow0 \quad \mathrm{a.e.}\ \mathrm{in}\
Q_T,
\end{cases}\ee as $\delta\rightarrow0^+$ (taking a
subsequence if necessary), {where $d_\delta=\rho_\delta+n_\delta$. (\ref{ndconvergence}) and Egrov theorem
imply that for any small positive constant $\sigma$, there exists
a domain $Q_T^\prime\subset Q_T$, such that
$|Q_T/Q_T^\prime|\le\sigma$ and that \be\label{3-a4}\begin{cases}
n_\delta-Ad_\delta\rightarrow0 \quad
\mathrm{uniformly}\ \mathrm{in}\ Q_T^\prime,\\[2mm]
\rho_\delta-Bd_\delta\rightarrow0 \quad \mathrm{uniformly}\
\mathrm{in}\ Q_T^\prime
\end{cases}\ee as
$\delta\rightarrow0^+$ $\Big($taking the same sequence as in
(\ref{ndconvergence})$\Big)$.

In view of (\ref{3-a4}), we obtain that there exists a positive
constant $\delta_0$ such that \be\label{ddeltanrho}
\begin{cases}
Ad_\delta\le n_\delta+1,\\[2mm]
Bd_\delta\le \rho_\delta+1
\end{cases}
\ee
 for $\delta\le \delta_0$ and any $(x,t)\in Q_T^\prime$. Note that
 $\delta_0$ does not depend on $(x,t)$.

Therefore for $\delta\le \delta_0$, $Ad_\delta$ and $Bd_\delta$
are bounded in $L^{\Gamma+\theta_1}(Q_T^\prime)$ and in
$L^{\gamma+\theta_2}(Q_T^\prime)$, respectively. Note that when
$\Gamma+\theta_1>\gamma+\theta_2$
or $\Gamma+\theta_1<\gamma+\theta_2$, one can not generally
guarantee that \bex\begin{cases}d_\delta=\rho_\delta+n_\delta\in
L^{\Gamma+\theta_1}(Q_T),\\[2mm]d_\delta=\rho_\delta+n_\delta\in
L^{\gamma+\theta_2}(Q_T),\end{cases}\eex since the only useful information
we have is \bex\begin{cases}
\rho_\delta\in L^{\gamma+\theta_2}(Q_T)\cap
L^\infty\Big(0,T;L^\gamma(\Omega)\Big),\\[2mm]n_\delta\in
L^{\Gamma+\theta_1}(Q_T)\cap
L^\infty\Big(0,T;L^\Gamma(\Omega)\Big).\end{cases}\eex
 Thus it indicates
that the weighted functions $A$ and $B$ can cancel some possible
oscillation of $d_\delta$.

Without loss of generality, we only show the proof of
(\ref{3-c1})$_1$. In fact, the proof of (\ref{3-c1})$_2$ is
similar. To begin with, we divide an integral into a sum of two
parts, i.e., {\it Integrability Part} + {\it Small Region Part}.
More precisely, we have \bex
\begin{split}\int_{Q_T}\Phi
T_k(\rho_\delta)P(n_\delta,\rho_\delta)\,dx\,dt
=\int_{Q^\prime_T}\Phi
T_k(\rho_\delta)P(n_\delta,\rho_\delta)\,dx\,dt +
\int_{Q_T/Q^\prime_T}\Phi
T_k(\rho_\delta)P(n_\delta,\rho_\delta)\,dx\,dt,
\end{split}
\eex for any $\Phi\in C(\overline{Q_T})$ where $\Phi\ge0$.

\subsubsection*{{\bf $\bullet$ Analysis of the  {\it Integrability Part}.}}

\be\label{III123}\begin{split}
&\lim\limits_{\delta\rightarrow0^+}\int_{Q^\prime_T}\Phi
 T_k(\rho_\delta) P(n_\delta,\rho_\delta)\,dx\,dt\\
=&\lim\limits_{\delta\rightarrow0^+}\int_{Q^\prime_T}\Phi
T_k(Bd_\delta)P(Ad_\delta,Bd_\delta)\,dx\,dt\\&+\lim\limits_{\delta\rightarrow0^+}\int_{Q^\prime_T}\Phi
\big[T_k(\rho_\delta)-T_k(Bd_\delta)\big]P(Ad_\delta,Bd_\delta)\,dx\,dt\\&+\lim\limits_{\delta\rightarrow0^+}\int_{Q^\prime_T}\Phi
T_k(\rho_\delta)\big[P(n_\delta,\rho_\delta)-P(Ad_\delta,Bd_\delta)\big]\,dx\,dt
\\=& \sum\limits_{i=1}^3III_i.
\end{split}
\ee

For $III_2$, in view of (\ref{3-a4}), the continuity of the map
$z\mapsto T_k(z)$, and the boundedness of $P(Ad_\delta,Bd_\delta)$
in $L^{\theta_m}(Q_T)$ due to (\ref{pressure ineq}),
(\ref{ddeltanrho}), and (\ref{3-hierho}), we have \be\label{III2}
III_2\rightarrow0 \ee as $\delta\rightarrow0^+$, where
$\theta_m=\min\{\frac{\Gamma+\theta_1}{\Gamma},
\frac{\gamma+\theta_2}{\gamma}\}$.

For $III_3$, similar to (\ref{pressurediff}) and
(\ref{rho+bound}), we get \be\label{pressurediff+1}
\begin{split}
&\Big|P(A_\delta d_\delta, B_\delta d_\delta)-P(A d_\delta,B
d_\delta)\Big|\\
\le&\frac{A_+\gamma\big(\frac{A_-}{A_+}\big)^\frac{1}{\Gamma}[\rho_+(\eta_1,\eta_2)]^{\gamma-\frac{\gamma}{\Gamma}}}{\frac{\gamma}{\Gamma}
[1-\alpha(\eta_1,\eta_2)]
 +\alpha(\eta_1,\eta_2)}\big|A_\delta
d_\delta-Ad_\delta\big|+\frac{A_+\gamma[\rho_+(\eta_1,\eta_2)]^{\gamma-1}}{\frac{\gamma}{\Gamma}[1-\alpha(\eta_1,\eta_2)]
 +\alpha(\eta_1,\eta_2)}\big|B_\delta
d_\delta-Bd_\delta\big|\\
\le&C[\rho_+(\eta_1,\eta_2)]^{\gamma-\frac{\gamma}{\Gamma}}\big|A_\delta
d_\delta-Ad_\delta\big|+
C[\rho_+(\eta_1,\eta_2)]^{\gamma-1}\big|B_\delta
d_\delta-Bd_\delta\big|\\
\le&C\Big[(A_\delta d_\delta+Ad_\delta)^{\Gamma-1}+(B_\delta
d_\delta+ B
d_\delta)^{\gamma(1-\frac{1}{\Gamma})}\Big]\big|A_\delta
d_\delta-Ad_\delta\big|\\&+ C\Big[(A_\delta d_\delta+A
d_\delta)^{\Gamma(1-\frac{1}{\gamma})}+(B_\delta d_\delta+ B
d_\delta)^{\gamma-1}\Big]\big|B_\delta d_\delta-Bd_\delta\big|,
\end{split}
\ee where we have used
 \be\label{rho+bound+1}
 \begin{cases}
\rho_+(\eta_1,\eta_2)\le
C_0^\frac{1}{\gamma}\Big(\eta_1^\frac{\Gamma}{\gamma}+\eta_2\Big),\\[2mm]
\eta_1\le A_\delta d_\delta+A d_\delta,\\[2mm]
\eta_2\le B_\delta d_\delta+ B d_\delta.
 \end{cases}
 \ee
Therefore we obtain \be\label{III3}
\begin{split}
|III_3|\le& C_k\lim\limits_{\delta\rightarrow0^+}\int_{Q^\prime_T}
\Big[(A_\delta d_\delta+Ad_\delta)^{\Gamma-1}+(B_\delta d_\delta+
B
d_\delta)^{\gamma(1-\frac{1}{\Gamma})}\Big]\big|n_\delta-Ad_\delta\big|\,dx\,dt\\+&C_k\lim\limits_{\delta\rightarrow0^+}\int_{Q^\prime_T}\Big[(A_\delta
d_\delta+A d_\delta)^{\Gamma(1-\frac{1}{\gamma})}+(B_\delta
d_\delta+ B
d_\delta)^{\gamma-1}\Big]\big|\rho_\delta-Bd_\delta\big|\,dx\,dt\\
&\rightarrow0
\end{split}
\ee as $\delta\rightarrow0^+$, due to (\ref{3-hierho}),
(\ref{3-a4}), and (\ref{ddeltanrho}).

In view of (\ref{III2}) and (\ref{III3}), (\ref{III123}) can be
refined as follows. \be\label{III123+1}\begin{split}
\lim\limits_{\delta\rightarrow0^+}\int_{Q^\prime_T}\Phi
 T_k(\rho_\delta)P(n_\delta,\rho_\delta)\,dx\,dt
=&\int_{Q^\prime_T}\Phi
\overline{\overline{T_k(Bd)P(Ad,Bd)}}\,dx\,dt\\
\ge&\int_{Q^\prime_T}\Phi
\overline{\overline{T_k(Bd)}}\,\,\overline{\overline{P(Ad,Bd)}}\,dx\,dt
\end{split}
\ee due to Lemma \ref{Feireisl-Novotny lemma} and the fact that
the maps $z\mapsto T_k(Bz)$ and $z\mapsto P(Az,Bz)$ are non
decreasing.

Note that \be\label{3-a10}\begin{split} &\int_{Q^\prime_T}\Phi
\overline{\overline{T_k(Bd)}}\,\,\overline{\overline{P(Ad,Bd)}}\,dx\,dt\\
=&\lim\limits_{\delta\rightarrow0^+}\int_{Q^\prime_T}\Phi
T_k(Bd_\delta)\,\,\overline{\overline{P(Ad,Bd)}}\,dx\,dt\\=&
\lim\limits_{\delta\rightarrow0^+}\int_{Q^\prime_T}\Phi
T_k(\rho_\delta)\,\,\overline{\overline{P(Ad,Bd)}}\,dx\,dt\\&+\lim\limits_{\delta\rightarrow0^+}\int_{Q^\prime_T}\Phi
\big[T_k(Bd_\delta)-T_k(\rho_\delta)\big]\,\,\overline{\overline{P(Ad,Bd)}}\,dx\,dt\\=&
\int_{Q^\prime_T}\Phi
\overline{T_k(\rho)}\,\,\overline{\overline{P(Ad,Bd)}}\,dx\,dt,
\end{split}
\ee where we have used (\ref{3-a4}), the continuity of the map
$z\mapsto T_k(z)$, and $\overline{\overline{P(Ad,Bd)}}\in
L^{\theta_m}(Q_T)$ with
$\theta_m=\min\{\frac{\Gamma+\theta_1}{\Gamma},
\frac{\gamma+\theta_2}{\gamma}\}>1$, such that \bex
\lim\limits_{\delta\rightarrow0^+}\int_{Q^\prime_T}\Phi
\big[T_k(Bd_\delta)-T_k(\rho_\delta)\big]\,\,\overline{\overline{P(Ad,Bd)}}\,dx\,dt\rightarrow0
\eex as $\delta\rightarrow0^+$. Similarly, we have
\be\label{3-a9}\begin{split} \int_{Q^\prime_T}\Phi
\overline{T_k(\rho)}\,\,\overline{\overline{P(Ad,Bd)}}\,dx\,dt
=&\lim\limits_{\delta\rightarrow0^+}\int_{Q^\prime_T}\Phi
\overline{T_k(\rho)}P(Ad_\delta,Bd_\delta)\,dx\,dt\\=&\lim\limits_{\delta\rightarrow0^+}\int_{Q^\prime_T}\Phi
\overline{T_k(\rho)}P(n_\delta,\rho_\delta)\,dx\,dt\\&+\lim\limits_{\delta\rightarrow0^+}\int_{Q^\prime_T}\Phi
\overline{T_k(\rho)}\Big[P(Ad_\delta,Bd_\delta) -
P(n_\delta,\rho_\delta)\Big]\,dx\,dt\\=&\int_{Q^\prime_T}\Phi
\overline{T_k(\rho)}\,\,\overline{P(n,\rho)}\,dx\,dt.
\end{split}\ee
Combining (\ref{3-a10}) and (\ref{3-a9}), we have
\be\label{3-a11}\begin{split} &\int_{Q^\prime_T}\Phi
\overline{\overline{\big[T_k(Bd)\big]}}\,\,\overline{\overline{P(Ad,Bd)}}\,dx\,dt
=\int_{Q^\prime_T}\Phi
\overline{T_k(\rho)}\,\,\overline{P(n,\rho)}\,dx\,dt.
\end{split}
\ee

Note that the left term of (\ref{3-a11}) is exactly the same as
the right term of (\ref{III123+1}). Hence we obtain from
(\ref{III123+1}) and (\ref{3-a11}) that
 \be\label{3-a12}\begin{split}
\lim\limits_{\delta\rightarrow0^+}\int_{Q^\prime_T}\Phi
T_k(\rho_\delta)P(n_\delta,\rho_\delta)\,dx\,dt
\ge\int_{Q^\prime_T}\Phi
\overline{T_k(\rho)}\,\,\overline{P(n,\rho)}\,dx\,dt.
\end{split}
\ee

\subsubsection*{{\bf $\bullet$ Analysis of the {\it Small Region Part}.}}

For fixed $k$, we have \be\label{Small region}
\lim\limits_{\delta\rightarrow 0^+}\int_{Q_T/Q^\prime_T}\Phi
T_k(\rho_\delta)P(n_\delta,\rho_\delta)\,dx\,dt=\int_{Q_T/Q^\prime_T}\Phi
\overline{T_k(\rho)P(n,\rho)}\,dx\,dt, \ee since
$T_k(\rho_\delta)P(n_\delta,\rho_\delta)$ is bounded in
$L^{\theta_m}(Q_T)$ uniformly for $\delta>0$, where
$\theta_m=\min\{\frac{\Gamma+\theta_1}{\Gamma},
\frac{\gamma+\theta_2}{\gamma}\}>1$.

\medskip

\subsubsection*{{\bf $\bullet$ Analysis of the whole Part.}}

By virtue of (\ref{3-a12}) and (\ref{Small region}), we have
\be\label{3-2}
\begin{split}&\int_{Q_T}\Phi
\overline{T_k(\rho)P(n,\rho)}\,dx\,dt
\\=&\lim\limits_{\delta\rightarrow
0^+}\int_{Q^\prime_T}\Phi
T_k(\rho_\delta)P(n_\delta,\rho_\delta)\,dx\,dt +
\lim\limits_{\delta\rightarrow 0^+}\int_{Q_T/Q^\prime_T}\Phi
T_k(\rho_\delta)P(n_\delta,\rho_\delta)\,dx\,dt\\
\ge&\int_{Q^\prime_T}\Phi
\overline{T_k(\rho)}\,\,\overline{P(n,\rho)}\,dx\,dt +
\int_{Q_T/Q^\prime_T}\Phi
\overline{T_k(\rho)P(n,\rho)}\,dx\,dt\\
=&\int_{Q_T}\Phi
\overline{T_k(\rho)}\,\,\overline{P(n,\rho)}\,dx\,dt-\int_{Q_T/Q^\prime_T}\Phi
\overline{T_k(\rho)}\,\,\overline{P(n,\rho)}\,dx\,dt
\\&+ \int_{Q_T/Q^\prime_T}\Phi
\overline{T_k(\rho)P(n,\rho)}\,dx\,dt.
\end{split}
\ee

Since $|Q_T/Q_T^\prime|\le\sigma$, letting $\sigma$ go to zero, we
obtain that the last two terms on the right hand side of
(\ref{3-2}) will vanish. Hence we have \bex
\begin{split}\int_{Q_T}\Phi
\overline{T_k(\rho)P(n,\rho)}\,dx\,dt \ge\int_{Q_T}\Phi
\overline{T_k(\rho)}\,\,\overline{P(n,\rho)}\,dx\,dt.
\end{split}
\eex

Since $\Phi$ is arbitrary, we get (\ref{3-c1})$_1$. By using the arguments similar to the proof of (\ref{3-c1})$_1$, we get (\ref{3-c1})$_2$. Therefore we complete the
proof of the lemma.

\endpf

\begin{lemma}\label{3-le:3.7} Let
$(\rho_\delta,n_\delta,u_\delta)$ be the solution stated in
Proposition \ref{2-le:aweak solution} and $(\rho,n,u)$ be the
limit, then \be\label{3-8}
\lim\limits_{\delta\rightarrow0^+}\int_0^T\psi\int_\Omega\phi
H_\delta
\big[T_k(\rho_\delta)+T_k(n_\delta)\big]\,dx\,dt=\int_0^T\psi\int_\Omega\phi
\overline{H}\ \big[\overline{T_k(\rho)}+
\overline{T_k(n)}\big]\,dx\,dt, \ee for any $\psi\in
C_0^\infty(0,T)$ and $\phi\in
 C_0^\infty(\Omega)$, where
 \be\label{Hdelta}\begin{cases}
H_\delta:=P(n_\delta,\rho_\delta)-(2\mu+\lambda)\mathrm{div}u_\delta,
\\[2mm]
\overline{H}:=\overline{P(n,\rho)}-(2\mu+\lambda)\mathrm{div}u.
\end{cases}
\ee
\end{lemma}
\begin{remark}
Lemma \ref{3-le:3.7} is motivated by \cite{Feireisl,Lions}. The
statement of the lemma for the two-fluid model can be found in
\cite{Vasseur-Wen-Yu}.
\end{remark}

\bigskip

 To show the strong convergence of
$\rho_{\delta}$ and $n_{\delta}$, motivated by
\cite{Feireisl,Lions} (see also \cite{Vasseur-Wen-Yu}), we define
\bex
L_k(z)=\left\{\begin{array}{l} z\log z,\quad 0\le z\le k, \\
[3mm] z\log k+z\displaystyle\int_k^z\frac{T_k(s)}{s^2}\,ds, \quad
z\ge k,
\end{array}
\right. \eex satisfying \bex\begin{split} L_k(z)=\beta_kz-2k\ for\
all\ z\ge 3k,
\end{split}
\eex   where
$$
\beta_k=\log
k+\displaystyle\int_k^{3k}\frac{T_k(s)}{s^2}\,ds+\frac{2}{3}.
$$
We denote $b_k(z):=L_k(z)-\beta_kz$ where $b^\prime_k(z)=0$ for
all large $z$, and \be\label{3-id} b^\prime_k(z)z-b_k(z)=T_k(z).
\ee

Note that $\rho_\delta, n_\delta \in L^{2}(Q_T)$, $\rho,n\in
L^{2}(Q_T)$, and $u_\delta, u\in L^2(0,T;H_0^1(\Omega))$. Then using the same arguments as in \cite{Vasseur-Wen-Yu} where
Lemma \ref{3-le:3.7} is used, we arrive at \bex\begin{split}
&\int_\Omega
[\overline{L_k(\rho)}-L_k(\rho)+\overline{L_k(n)}-L_k(n)]\,dx
 \\
=&\frac{1}{2\mu+\lambda}\int_0^t\int_\Omega
\big(\overline{T_k(\rho)}
+\overline{T_k(n)}\big)\overline{P(n,\rho)}\,dx\,ds
\\&-\frac{1}{2\mu+\lambda}\int_0^t\int_\Omega
\overline{\big[T_k(\rho)+T_k(n)\big]P(n,\rho)}\,dx\,ds
\\&+\int_0^t\int_\Omega
[T_k(\rho)-\overline{T_k(\rho)}+T_k(n)-\overline{T_k(n)}]\mathrm{div}u\,dx\,ds.
\end{split}
\eex  This together with (\ref{3-c1}) yields
\be\label{3-last4}\begin{split} \int_\Omega
[\overline{L_k(\rho)}-L_k(\rho)+\overline{L_k(n)}-L_k(n)]\,dx
\le\int_0^t\int_\Omega
[T_k(\rho)-\overline{T_k(\rho)}+T_k(n)-\overline{T_k(n)}]\mathrm{div}u\,dx\,ds.
\end{split}
\ee

\bigskip

In order to include the case that both $\gamma$ and
$\Gamma$ can touch $\frac{9}{5}$, we need the following new estimate.
\begin{lemma}\label{3-le5.7}Let
$(\rho_\delta,n_\delta)$ be the solution stated in Proposition
\ref{2-le:aweak solution} and $(\rho,n)$ be the limit, then \bex
\lim\limits_{\delta\rightarrow0}\|T_k(Ad_\delta) + T_k(Bd_\delta)
-
T_k(Ad)-T_k(Bd)\|^{\Gamma_{min}+1}_{L^{\Gamma_{min}+1}(Q^\prime_T)}\le
C_k\sigma^\frac{K_{min}-1}{K_{min}}+C\eex for any
$\Gamma,\gamma\ge\frac{9}{5}$, and any given $k>0$, where $C$ is
independent of $\sigma$, $\delta$, and $k$, and $C_k$ is
independent of $\sigma$ and $\delta$ but may depend on $k$. Here
\be\label{3-Gamma-K-min}\Gamma_{min}=\min\{\Gamma,\gamma\},\quad
K_{min}=\min\{\frac{\Gamma+\theta_1}{\Gamma},
\frac{\gamma+\theta_2}{\gamma},2\}. \ee
\end{lemma}

 \pf Note that \be\label{3-Tk-Tk}\begin{split}&
|T_k(Ad_\delta) + T_k(Bd_\delta) -
T_k(Ad)-T_k(Bd)|^{\Gamma_{min}+1}\\ \le&
\Big(A|d_\delta-d|+B|d_\delta-d|\Big)^{\Gamma_{min}}|T_k(Ad_\delta)
+ T_k(Bd_\delta) - T_k(Ad)-T_k(Bd)|\\ \le&
C\Big(|A(d_\delta-d)|^{\Gamma} +
|B(d_\delta-d)|^{\gamma}+1\Big)\big|T_k(Ad_\delta) +
T_k(Bd_\delta) - T_k(Ad)-T_k(Bd)\big| \\ \le&
C\Big[(Ad_\delta)^\Gamma-(Ad)^\Gamma +
(Bd_\delta)^\gamma-(Bd)^{\gamma}\Big]\big[T_k(Ad_\delta) +
T_k(Bd_\delta) - T_k(Ad)-T_k(Bd)\big]\\&+C\big|T_k(Ad_\delta) +
T_k(Bd_\delta) - T_k(Ad)-T_k(Bd)\big| \end{split}\ee due to the
fact that \bex |T^\prime(x)|\le 1 \eex for any $x\ge0$, and that
\bex |x-y|^\Gamma\le |x^\Gamma-y^\Gamma|,\quad |x-y|^\gamma\le
|x^\gamma-y^\gamma|
 \eex for any $x,y\ge0$.

Therefore we have \be\label{3-total}\begin{split}
&\int_{Q^\prime_T}|T_k(Ad_\delta) + T_k(Bd_\delta) -
T_k(Ad)-T_k(Bd)|^{\Gamma_{min}+1}\,dx\,dt
\\ \le&C\int_{Q^\prime_T}\Big[(Ad_\delta)^\Gamma-(Ad)^\Gamma +
(Bd_\delta)^\gamma-(Bd)^{\gamma}\Big]\big[T_k(Ad_\delta) +
T_k(Bd_\delta) -
T_k(Ad)-T_k(Bd)\big]\,dx\,dt\\&+C\int_{Q^\prime_T}\big|T_k(Ad_\delta)
+ T_k(Bd_\delta) - T_k(Ad)-T_k(Bd)\big|\,dx\,dt\\:=&
IV^\delta_1+IV^\delta_2,
\end{split}
\ee where $Q^\prime_T$ is introduced in (\ref{3-a4}).

For $IV_1^\delta$, we have \bex\begin{split}
IV^\delta_1=&C\int_{Q^\prime_T}\Big[(Ad_\delta)^\Gamma
+(Bd_\delta)^\gamma\Big]\big[T_k(Ad_\delta) +
T_k(Bd_\delta)\big]\,dx\,dt\\&-C\int_{Q^\prime_T}\Big[(Ad_\delta)^\Gamma
+(Bd_\delta)^\gamma\Big]\big[T_k(Ad)+T_k(Bd)\big]\,dx\,dt\\&-C\int_{Q^\prime_T}\Big[(Ad)^\Gamma+(Bd)^{\gamma}\Big]\big[T_k(Ad_\delta)
+
T_k(Bd_\delta)\big]\,dx\,dt\\&+C\int_{Q^\prime_T}\Big[(Ad)^\Gamma+(Bd)^{\gamma}\Big]\big[T_k(Ad)+T_k(Bd)\big]\,dx\,dt.
\end{split}
\eex
Taking the limit as $\delta\rightarrow0$ (taking a
subsequence if necessary), we have
 \bex \begin{split}
\lim\limits_{\delta\rightarrow0}IV^\delta_1
=&C\lim\limits_{\delta\rightarrow0}\int_{Q^\prime_T}\Big[(Ad_\delta)^\Gamma
+(Bd_\delta)^\gamma\Big]\big[T_k(Ad_\delta) +
T_k(Bd_\delta)\big]\,dx\,dt\\&-C\int_{Q^\prime_T}\overline{\overline{(Ad)^\Gamma
+(Bd)^\gamma}}\,\,\overline{\overline{T_k(Ad) +
T_k(Bd)}}\,dx\,dt\\&+C\int_{Q^\prime_T}\big[\overline{\overline{(Ad)^\Gamma
+(Bd)^\gamma}}-(Ad)^\Gamma-(Bd)^{\gamma}\big]\big[\overline{\overline{T_k(Ad)
+ T_k(Bd)}}-T_k(Ad)-T_k(Bd)\big]\,dx\,dt.
\end{split}
\eex Due to the convexity of $z\mapsto (Bz)^\gamma+(Az)^\Gamma$
and the concavity of $z\mapsto T_k(Az) + T_k(Bz)$ such that
\bex\begin{cases} \overline{\overline{(Ad)^\Gamma
+(Bd)^\gamma}}\ge (Ad)^\Gamma+(Bd)^{\gamma},\\[2mm]
\overline{\overline{T_k(Ad) + T_k(Bd)}}\le T_k(Ad) + T_k(Bd),
\end{cases}
\eex we have  \be\label{3-IV1} \begin{split}
\lim\limits_{\delta\rightarrow0}IV^\delta_1
\le&C\lim\limits_{\delta\rightarrow0}\int_{Q^\prime_T}\Big[(Ad_\delta)^\Gamma
+(Bd_\delta)^\gamma\Big]\big[T_k(Ad_\delta) +
T_k(Bd_\delta)\big]\,dx\,dt\\&-C\int_{Q^\prime_T}\overline{\overline{(Ad)^\Gamma
+(Bd)^\gamma}}\,\,\overline{\overline{T_k(Ad) + T_k(Bd)}}\,dx\,dt.
\end{split}
\ee

For $IV^\delta_2$, we apply Young inequality and obtain
\be\label{3-IV2}
IV^\delta_2\le\frac{1}{2}\int_{Q^\prime_T}|T_k(Ad_\delta) +
T_k(Bd_\delta) - T_k(Ad)-T_k(Bd)|^{\Gamma_{min}+1}\,dx\,dt +C_1.
\ee

Combining (\ref{3-total}) with (\ref{3-IV1}) and (\ref{3-IV2})
yields \be\label{3-total1}\begin{split}
&\lim\limits_{\delta\rightarrow0}\int_{Q^\prime_T}|T_k(Ad_\delta)
+ T_k(Bd_\delta) - T_k(Ad)-T_k(Bd)|^{\Gamma_{min}+1}\,dx\,dt
\\ \le&C\lim\limits_{\delta\rightarrow0}\int_{Q^\prime_T}\Big[(Ad_\delta)^\Gamma
+(Bd_\delta)^\gamma\Big]\big[T_k(Ad_\delta) +
T_k(Bd_\delta)\big]\,dx\,dt\\&-C\int_{Q^\prime_T}\overline{\overline{(Ad)^\Gamma
+(Bd)^\gamma}}\,\,\overline{\overline{T_k(Ad) +
T_k(Bd)}}\,dx\,dt+C_1.
\end{split}
\ee
On the other hand, \be\label{3-a13}
\begin{split}
\frac{\partial P(Az,Bz)}{\partial
z}=A_+\gamma\rho_+^{\gamma-1}(Az,Bz)\Big[\frac{\partial
\rho_+(Az,Bz)}{\partial n} A + \frac{\partial
\rho_+(Az,Bz)}{\partial \rho} B\Big].
\end{split}
\ee Recalling (\ref{partialrhon}), we have \bex
\begin{cases} \frac{\partial\rho_+(Az,Bz)}{\partial
n}=\frac{\big(\frac{A_-}{A_+}\big)^\frac{1}{\Gamma}[\rho_+(Az,Bz)]^{1-\frac{\gamma}{\Gamma}}}{\frac{\gamma}{\Gamma}(1-\alpha)
 +\alpha},\\[2mm]
\frac{\partial\rho_+(Az,Bz)}{\partial\rho}=\frac{1}{\frac{\gamma}{\Gamma}(1-\alpha)
 +\alpha}.
\end{cases}
 \eex
This together with (\ref{3-a13}) gives
\be\label{3-in}\begin{split} \frac{\partial P(Az,Bz)}{\partial
z}=&A_+\gamma\rho_+^{\gamma-1}(Az,Bz)\Big[\frac{\big(\frac{A_-}{A_+}\big)^\frac{1}{\Gamma}[\rho_+(Az,Bz)]^{1-\frac{\gamma}{\Gamma}}}{\frac{\gamma}{\Gamma}(1-\alpha)
 +\alpha} A + \frac{1}{\frac{\gamma}{\Gamma}(1-\alpha)
 +\alpha} B\Big]\\=&\frac{A_+\gamma}{\frac{\gamma}{\Gamma}(1-\alpha)
 +\alpha}\Big[\big(\frac{A_-}{A_+}\big)^\frac{1}{\Gamma}[\rho_+(Az,Bz)]^{\gamma-\frac{\gamma}{\Gamma}}A  + \rho_+^{\gamma-1}(Az,Bz)B\Big]
 \\ \ge&\frac{A_+\gamma}{\max\{\frac{\gamma}{\Gamma},1\}}\Big[\big(\frac{A_-}{A_+}\big)^\frac{1}{\Gamma}[{\rho^\gamma}_+(Az,Bz)]^{1-\frac{1}{\Gamma}}A
+ [\rho_+^{\gamma}(Az,Bz)]^{1-\frac{1}{\gamma}}B\Big]\\
\ge&C_2\Big[A^\Gamma z^{\Gamma-1} + B^\gamma z^{\gamma-1}\Big],
\end{split}
\ee where we have used (\ref{pressure ineq}), and
$C_2=C_2(A_+,A_-,\Gamma,\gamma)>0$.

Thus, we introduce \be\label{3-increasing function}\begin{split}
G_{A,B}(z):=&P(Az,Bz)-\frac{C_2}{\max\{\Gamma,\gamma\}}\Big[(Az)^{\Gamma}
+ (Bz)^{\gamma}\Big],
\end{split}
\ee which is inspired by \cite{FeireislJDE} for the single-phase
flow where non-mono pressure of one component is studied.

In view of (\ref{3-in}) and (\ref{3-increasing function}), we
obtain \bex \frac{d}{dz}G_{A,B}(z)= \frac{\partial
P(Az,Bz)}{\partial
z}-C_2\Big[\frac{\Gamma}{\max\{\Gamma,\gamma\}}A^{\Gamma}z^{\Gamma-1}
+
\frac{\gamma}{\max\{\Gamma,\gamma\}}B^{\gamma}z^{\gamma-1}\Big]\ge
0,\eex and thus $z\mapsto G_{A,B}(z)$ is a non-decreasing function
over $[0,\infty)$.

Let's return to (\ref{3-total1}), and make use of $G_{A,B}(z)$. Then
we get \be\label{3-total2}\begin{split}
&\lim\limits_{\delta\rightarrow0}\int_{Q^\prime_T}|T_k(Ad_\delta)
+ T_k(Bd_\delta) - T_k(Ad)-T_k(Bd)|^{\Gamma_{min}+1}\,dx\,dt
\\ \le&\frac{C\max\{\Gamma,\gamma\}}{C_2}\lim\limits_{\delta\rightarrow0}\int_{Q^\prime_T}\frac{C_2}{\max\{\Gamma,\gamma\}}\Big[(Ad_\delta)^\Gamma
+(Bd_\delta)^\gamma\Big]\big[T_k(Ad_\delta) +
T_k(Bd_\delta)\big]\,dx\,dt\\&-\frac{C\max\{\Gamma,\gamma\}}{C_2}\int_{Q^\prime_T}\overline{\overline{\frac{C_2}{\max\{\Gamma,\gamma\}}\big[(Ad)^\Gamma
+(Bd)^\gamma\big]}}\,\,\overline{\overline{T_k(Ad) +
T_k(Bd)}}\,dx\,dt+C_1\\
=&\frac{C\max\{\Gamma,\gamma\}}{C_2}\lim\limits_{\delta\rightarrow0}\int_{Q^\prime_T}P(Ad_\delta,Bd_\delta)\big[T_k(Ad_\delta)
+
T_k(Bd_\delta)\big]\,dx\,dt\\&-\frac{C\max\{\Gamma,\gamma\}}{C_2}\int_{Q^\prime_T}\overline{\overline{P(Ad,Bd)}}\,\,\overline{\overline{T_k(Ad)
+
T_k(Bd)}}\,dx\,dt\\&-\frac{C\max\{\Gamma,\gamma\}}{C_2}\lim\limits_{\delta\rightarrow0}\int_{Q^\prime_T}G_{A,B}(d_\delta)\big[T_k(Ad_\delta)
+
T_k(Bd_\delta)\big]\,dx\,dt\\&+\frac{C\max\{\Gamma,\gamma\}}{C_2}\int_{Q^\prime_T}\overline{\overline{G_{A,B}(d)}}\,\,\overline{\overline{T_k(Ad)+
T_k(Bd)}}\,dx\,dt+C_1.
\end{split}
\ee Note that \be\label{3-5.38}\begin{split}
&-\lim\limits_{\delta\rightarrow0}\int_{Q^\prime_T}G_{A,B}(d_\delta)\big[T_k(Ad_\delta)
+
T_k(Bd_\delta)\big]\,dx\,dt+\int_{Q^\prime_T}\overline{\overline{G_{A,B}(d)}}\,\,\overline{\overline{T_k(Ad)+
T_k(Bd)}}\,dx\,dt\\&=\int_{Q^\prime_T}\overline{\overline{G_{A,B}(d)}}\,\,\overline{\overline{T_k(Ad)+
T_k(Bd)}}\,dx\,dt-\int_{Q^\prime_T}\overline{\overline{G_{A,B}(d)\big[T_k(Ad)
+ T_k(Bd)\big]}}\,dx\,dt\\ &\le0,
\end{split}
\ee due to Lemma \ref{Feireisl-Novotny lemma} and the fact that
$z\mapsto G_{A,B}(z)$ and $z\mapsto T_k(Az) + T_k(Bz)$ are
non-decreasing functions.

By virtue of (\ref{3-5.38}), (\ref{3-total2}) yields
\be\label{3-total3}\begin{split}
&\lim\limits_{\delta\rightarrow0}\int_{Q^\prime_T}|T_k(Ad_\delta)
+ T_k(Bd_\delta) - T_k(Ad)-T_k(Bd)|^{\Gamma_{min}+1}\,dx\,dt
\\ \le& \frac{C\max\{\Gamma,\gamma\}}{C_2}\lim\limits_{\delta\rightarrow0}\int_{Q^\prime_T}P(Ad_\delta,Bd_\delta)\big[T_k(Ad_\delta)
+
T_k(Bd_\delta)\big]\,dx\,dt\\&-\frac{C\max\{\Gamma,\gamma\}}{C_2}\int_{Q^\prime_T}\overline{\overline{P(Ad,Bd)}}\,\,\overline{\overline{T_k(Ad)
+ T_k(Bd)}}\,dx\,dt+C_1.
\end{split}
\ee

By virtue of the uniform convergence (\ref{3-a4}), we rewrite
(\ref{3-total3}) as \be\label{3-total4}\begin{split}
&\lim\limits_{\delta\rightarrow0}\int_{Q^\prime_T}|T_k(Ad_\delta)
+ T_k(Bd_\delta) - T_k(Ad)-T_k(Bd)|^{\Gamma_{min}+1}\,dx\,dt
\\ \le& \frac{C\max\{\Gamma,\gamma\}}{C_2}\lim\limits_{\delta\rightarrow0}\int_{Q^\prime_T}P(n_\delta,\rho_\delta)\big[T_k(n_\delta)
+
T_k(\rho_\delta)\big]\,dx\,dt\\&-\frac{C\max\{\Gamma,\gamma\}}{C_2}\lim\limits_{\delta\rightarrow0}\int_{Q^\prime_T}P(Ad_\delta,Bd_\delta)\,\,\overline{\overline{T_k(Ad)
+ T_k(Bd)}}\,dx\,dt+C_1\\ =&
\frac{C\max\{\Gamma,\gamma\}}{C_2}\lim\limits_{\delta\rightarrow0}\int_{Q^\prime_T}P(n_\delta,\rho_\delta)\big[T_k(n_\delta)
+
T_k(\rho_\delta)\big]\,dx\,dt\\&-\frac{C\max\{\Gamma,\gamma\}}{C_2}\lim\limits_{\delta\rightarrow0}\int_{Q^\prime_T}P(n_\delta,\rho_\delta)\,\,\overline{\overline{T_k(Ad)
+ T_k(Bd)}}\,dx\,dt+C_1\\ =&
\frac{C\max\{\Gamma,\gamma\}}{C_2}\lim\limits_{\delta\rightarrow0}\int_{Q^\prime_T}P(n_\delta,\rho_\delta)\big[T_k(n_\delta)
+
T_k(\rho_\delta)\big]\,dx\,dt\\&-\frac{C\max\{\Gamma,\gamma\}}{C_2}\int_{Q^\prime_T}\overline{P(n,\rho)}\,\,\overline{\overline{T_k(Ad)
+ T_k(Bd)}}\,dx\,dt+C_1.
\end{split}
\ee

Similarly for the second term on the right hand side of
(\ref{3-total4}), we have \be\label{3-total5}\begin{split}
&\lim\limits_{\delta\rightarrow0}\int_{Q^\prime_T}|T_k(Ad_\delta)
+ T_k(Bd_\delta) - T_k(Ad)-T_k(Bd)|^{\Gamma_{min}+1}\,dx\,dt
 \\ \le&
\frac{C\max\{\Gamma,\gamma\}}{C_2}\lim\limits_{\delta\rightarrow0}\int_{Q^\prime_T}P(n_\delta,\rho_\delta)\big[T_k(n_\delta)
+
T_k(\rho_\delta)\big]\,dx\,dt\\&-\frac{C\max\{\Gamma,\gamma\}}{C_2}\int_{Q^\prime_T}\overline{P(n,\rho)}\,\,\overline{T_k(n)
+ T_k(\rho)}\,dx\,dt+C_1\\ =&
\frac{C\max\{\Gamma,\gamma\}}{C_2}\lim\limits_{\delta\rightarrow0}\int_{Q^\prime_T}H_\delta\big[T_k(n_\delta)+
T_k(\rho_\delta)\big]\,dx\,dt
\\&+\frac{C\max\{\Gamma,\gamma\}(2\mu+\lambda)}{C_2}\lim\limits_{\delta\rightarrow0}\int_{Q^\prime_T}\mathrm{div}u_\delta\big[T_k(n_\delta)+
T_k(\rho_\delta)\big]\,dx\,dt\\&-\frac{C\max\{\Gamma,\gamma\}}{C_2}\int_{Q^\prime_T}\overline{P(n,\rho)}\,\,\overline{T_k(n)
+ T_k(\rho)}\,dx\,dt+C_1.
\end{split}
\ee In view of (\ref{3-8}), we can take some appropriate test
functions, for example,
 \be\label{3-testfunction1}\psi_j\in
C_0^\infty(0,T),\quad \psi_j(t)\equiv 1\ \mathrm{for}\
\mathrm{any}\ t\in[\frac{1}{j},T-\frac{1}{j}],\ 0\le\psi_j\le1,\
\psi_j\rightarrow 1,\ee as $j\rightarrow\infty$, and
\be\label{3-testfunction2}\phi_j\in C_0^\infty(\Omega),\quad
\phi_j(x)\equiv1\ \mathrm{for}\ \mathrm{any}\ x\in\big\{x\in\Omega
\big| \mathrm{dist(x,\partial\Omega)\ge\frac{1}{j}}\big\},\
0\le\phi_j\le1,\ \phi_j\rightarrow 1,\ee as $j\rightarrow\infty$,
such that \be\label{3-a14}
\lim\limits_{\delta\rightarrow0^+}\int_{Q_T} H_\delta
\big[T_k(\rho_\delta)+T_k(n_\delta)\big]\,dx\,dt=\int_{Q_T}
\overline{H}\ \big[\overline{T_k(\rho)}+
\overline{T_k(n)}\big]\,dx\,dt. \ee Then from (\ref{3-total5}) and
(\ref{3-a14}), we obtain \be\label{3-total6}\begin{split}
&\lim\limits_{\delta\rightarrow0}\int_{Q^\prime_T}|T_k(Ad_\delta)
+ T_k(Bd_\delta) - T_k(Ad)-T_k(Bd)|^{\Gamma_{min}+1}\,dx\,dt
 \\ \le&
\frac{C\max\{\Gamma,\gamma\}}{C_2}\int_{Q_T}\overline{H}\,\,\overline{T_k(n)+
T_k(\rho)}\,dx\,dt-\frac{C\max\{\Gamma,\gamma\}}{C_2}\lim\limits_{\delta\rightarrow0}\int_{Q_T/Q^\prime_T}H_\delta\big[T_k(n_\delta)+
T_k(\rho_\delta)\big]\,dx\,dt
\\&+\frac{C\max\{\Gamma,\gamma\}(2\mu+\lambda)}{C_2}\lim\limits_{\delta\rightarrow0}\int_{Q^\prime_T}\mathrm{div}u_\delta\big[T_k(n_\delta)+
T_k(\rho_\delta)\big]\,dx\,dt\\&-\frac{C\max\{\Gamma,\gamma\}}{C_2}\int_{Q^\prime_T}\overline{P(n,\rho)}\,\,\overline{T_k(n)
+ T_k(\rho)}\,dx\,dt+C_1,
\end{split}
\ee and thus \be\label{3-total7}\begin{split}
&\lim\limits_{\delta\rightarrow0}\int_{Q^\prime_T}|T_k(Ad_\delta)
+ T_k(Bd_\delta) - T_k(Ad)-T_k(Bd)|^{\Gamma_{min}+1}\,dx\,dt
 \\ \le&
\frac{C\max\{\Gamma,\gamma\}}{C_2}\int_{Q_T/Q^\prime_T}\overline{P(n,\rho)}\,\,\overline{T_k(n)+
T_k(\rho)}\,dx\,dt-\frac{C\max\{\Gamma,\gamma\}}{C_2}\int_{Q_T/Q^\prime_T}\overline{H\big[T_k(n)+
T_k(\rho)\big]}\,dx\,dt
\\&+\frac{C\max\{\Gamma,\gamma\}(2\mu+\lambda)}{C_2}\lim\limits_{\delta\rightarrow0}\int_{Q^\prime_T}\mathrm{div}u_\delta\big[T_k(n_\delta)+
T_k(\rho_\delta)\big]\,dx\,dt\\&-\frac{C\max\{\Gamma,\gamma\}(2\mu+\lambda)}{C_2}\int_{Q_T}\mathrm{div}u\overline{T_k(n)+
T_k(\rho)}\,dx\,dt+C_1 \\ =&
\frac{C\max\{\Gamma,\gamma\}}{C_2}\int_{Q_T/Q^\prime_T}\Big[\overline{P(n,\rho)}\,\,\overline{T_k(n)+
T_k(\rho)}-\overline{H\big[T_k(n)+ T_k(\rho)\big]}\Big]\,dx\,dt
\\&+\frac{C\max\{\Gamma,\gamma\}(2\mu+\lambda)}{C_2}\lim\limits_{\delta\rightarrow0}\int_{Q^\prime_T}\mathrm{div}u_\delta\Big[T_k(n_\delta)+
T_k(\rho_\delta)-\overline{T_k(n)+
T_k(\rho)}\Big]\,dx\,dt+C_1\\:=&V_1+V_2+C_1,
\end{split}
\ee since \bex\begin{split}
H_\delta:=&P(n_\delta,\rho_\delta)-(2\mu+\lambda)\mathrm{div}u_\delta,
\\
\overline{H}:=&\overline{P(n,\rho)}-(2\mu+\lambda)\mathrm{div}u.
\end{split}
\eex

For $V_1$, we apply H\"older inequality, (\ref{3-lim})$_3$, and
(\ref{3-hierho}), and then obtain \be\label{3-V1}\begin{split}
V_1\le& \frac{C\max\{\Gamma,\gamma\}}{C_2}
\Big\{\int_{Q_T/Q^\prime_T}\Big|\overline{P(n,\rho)}\,\,\overline{T_k(n)+
T_k(\rho)}-\overline{H\big[T_k(n)+
T_k(\rho)\big]}\Big|^{K_{min}}\,dx\,dt\Big\}^\frac{1}{K_{min}}\Big|Q_T/Q^\prime_T\Big|^\frac{K_{min}-1}{K_{min}}\\
\le& C^k_3\sigma^\frac{K_{min}-1}{K_{min}},
\end{split}
\ee where $K_{min}=\min\{\frac{\Gamma+\theta_1}{\Gamma},
\frac{\gamma+\theta_2}{\gamma},2\}>1$, and $C^k_3$ is independent
of $\sigma$ for $\sigma\in(0,1)$ but may depend on $k$.

For $V_2$, by virtue of H\"older inequality and (\ref{2-r4}), we
have \bex\begin{split} V_2\le&
\frac{C\max\{\Gamma,\gamma\}(2\mu+\lambda)}{C_2}\limsup\limits_{\delta\rightarrow0}\Big(\int_{Q^\prime_T}|\mathrm{div}u_\delta|^2\,dx\,dt\Big)^\frac{1}{2}\Big(\int_{Q^\prime_T}\Big|T_k(n_\delta)+
T_k(\rho_\delta)-\overline{T_k(n)+
T_k(\rho)}\Big|^2\,dx\,dt\Big)^\frac{1}{2}\\ \le&
C_4\Big(\limsup\limits_{\delta\rightarrow0}\int_{Q^\prime_T}\Big|T_k(n_\delta)+
T_k(\rho_\delta)-T_k(n)-
T_k(\rho)\Big|^2\,dx\,dt\Big)^\frac{1}{2}\\
&+C_4\Big(\limsup\limits_{\delta\rightarrow0}\int_{Q^\prime_T}\Big|\overline{T_k(n)+
T_k(\rho)}-T_k(n)- T_k(\rho)\Big|^2\,dx\,dt\Big)^\frac{1}{2},
\end{split}
\eex where $C_4$ is independent of $\sigma$, $\delta$, and $k$.
This together with the lower semi-continuity of $L^2$ norm and
Young inequality deduces that \be\label{3-V2}\begin{split} V_2\le&
\frac{1}{2}\limsup\limits_{\delta\rightarrow0}\int_{Q^\prime_T}\Big|T_k(n_\delta)+
T_k(\rho_\delta)-T_k(n)- T_k(\rho)\Big|^{\Gamma_{min}+1}\,dx\,dt+C_5\\
=&
\frac{1}{2}\lim\limits_{\delta\rightarrow0}\int_{Q^\prime_T}\Big|T_k(Ad_\delta)+
T_k(Bd_\delta)-T_k(Ad)- T_k(Bd)\Big|^{\Gamma_{min}+1}\,dx\,dt+C_5,
\end{split}
\ee for some positive constant $C_5$ independent of $\sigma$,
$\delta$, and $k$. Here we have applied (\ref{ndconvergence}) to
the equality.

Thus we have \bex\begin{split}
\lim\limits_{\delta\rightarrow0}\int_{Q^\prime_T}|T_k(Ad_\delta) +
T_k(Bd_\delta) - T_k(Ad)-T_k(Bd)|^{\Gamma_{min}+1}\,dx\,dt
  \le
2C^k_3\sigma^\frac{K_{min}-1}{K_{min}}+2C_5+2C_1,
\end{split}
\eex according to (\ref{3-total7}), (\ref{3-V1}), and
(\ref{3-V2}).

The proof of the lemma is complete.

\endpf

\begin{corollary}\label{3-coromeasure} Let
$(\rho_\delta,n_\delta)$ be the solution stated in Proposition
\ref{2-le:aweak solution} and $(\rho,n)$ be the limit, then \bex
\lim\limits_{\delta\rightarrow0}\|T_k(n_\delta) + T_k(\rho_\delta)
- T_k(n)-T_k(\rho)\|_{L^{\Gamma_{min}+1}(Q_T)}^{\Gamma_{min}+1}\le
C\eex for any given $k>0$, where $C$ is independent of $\sigma$,
$\delta$, and $k$. Here $\Gamma_{min}$ and $K_{min}$ are given by
(\ref{3-Gamma-K-min}).
\end{corollary}
\pf In view of (\ref{3-a4}), we have
\bex\begin{split}&\lim\limits_{\delta\rightarrow0}\int_{Q_T}|T_k(n_\delta)
+ T_k(\rho_\delta) -
T_k(n)-T_k(\rho)|^{\Gamma_{min}+1}\,dx\,dt\\=&\lim\limits_{\delta\rightarrow0}\int_{Q_T/Q^\prime_T}|T_k(n_\delta)
+ T_k(\rho_\delta) -
T_k(n)-T_k(\rho)|^{\Gamma_{min}+1}\,dx\,dt\\&+
\lim\limits_{\delta\rightarrow0}\int_{Q^\prime_T}|T_k(Ad_\delta) +
T_k(Bd_\delta) - T_k(Ad)-T_k(Bd)|^{\Gamma_{min}+1}\,dx\,dt.
\end{split}
\eex Similar to (\ref{3-V1}), the first term on the right hand
side will tend to zero as $\sigma\rightarrow0^+$. And for the
second term, we use Lemma \ref{3-le5.7}. Consequently, letting
$\sigma\rightarrow0^+$, we complete the proof of the corollary.

\endpf

Corollary \ref{3-coromeasure} combined with the lower
semi-continuity of the norm yields the following corollary.
\begin{corollary}\label{3-cor5.7}Let
$(\rho_\delta,n_\delta)$ be the solution stated in Proposition
\ref{2-le:aweak solution} and $(\rho,n)$ be the limit, then \bex
\|\overline{T_k(n)} + \overline{T_k(\rho)} -
T_k(n)-T_k(\rho)\|_{L^{\Gamma_{min}+1}(Q_T)}\le C\eex for any
given $k>0$, where $C$ is independent of $k$.
\end{corollary}

\endpf

\bigskip

Denote \be\label{3-10} Q_{T,k}=\Big\{(x,t)\in Q_T\big|
\rho(x,t)\ge k,\, \mathrm{or}\,\, n(x,t)\ge k\Big\}. \ee

Here we are able  to control the right-hand side of
\eqref{3-last4} in the following lemma.
\begin{lemma}Let
$(\rho_\delta,n_\delta,u_\epsilon)$ be the solution stated in
Proposition \ref{2-le:aweak solution} and $(\rho,n,u)$ be the
limit, then \label{right side lemma}
 \be\label{3-c2} \begin{split}
 \lim\limits_{k\rightarrow\infty}\int_{Q_T}
[T_k(\rho)-\overline{T_k(\rho)}+T_k(n)-\overline{T_k(n)}]\mathrm{div}u\,dx\,dt=0.
\end{split}
\ee
\end{lemma}
\pf Using H\"older inequality and Corollary \ref{3-cor5.7}, we
have\be\label{3-a15}\begin{split} &\int_{Q_T}
[T_k(\rho)-\overline{T_k(\rho)}+T_k(n)-\overline{T_k(n)}]\mathrm{div}u\,dx\,dt\\
=&\int_{Q_{T,k}}
[T_k(\rho)-\overline{T_k(\rho)}+T_k(n)-\overline{T_k(n)}]\mathrm{div}u\,dx\,dt\\&+
\int_{Q_T/Q_{T,k}}
[T_k(\rho)-\overline{T_k(\rho)}+T_k(n)-\overline{T_k(n)}]\mathrm{div}u\,dx\,dt\\
\le&\|T_k(\rho)-\overline{T_k(\rho)}+T_k(n)-\overline{T_k(n)}\|_{L^2(Q_{T,k})}\|\mathrm{div}u\|_{L^2(Q_{T,k})}\\&+
\|T_k(\rho)-\overline{T_k(\rho)}+T_k(n)-\overline{T_k(n)}\|_{L^2(Q_T/Q_{T,k})}\|\mathrm{div}u\|_{L^2(Q_T/Q_{T,k})}\\
\le&C\|\mathrm{div}u\|_{L^2(Q_{T,k})}+C\|T_k(\rho)-\overline{T_k(\rho)}+T_k(n)-\overline{T_k(n)}\|_{L^2(Q_T/Q_{T,k})}.
\end{split}
\ee For the second term on the right hand side of (\ref{3-a15}),
by virtue of the standard interpolation inequality and Corollary
\ref{3-cor5.7}, we have \be\label{3-a16}\begin{split}&\|T_k(\rho)-
\overline{T_k(\rho)}+T_k(n)-
\overline{T_k(n)}\|_{L^2(Q_T/Q_{T,k})}\\
 \le& \|T_k(\rho)- \overline{T_k(\rho)}+T_k(n)-
\overline{T_k(n)}\|_{L^{1}(Q_T/Q_{T,k})}^{\frac{\Gamma_{min}-1}{2\Gamma_{min}}}\|T_k(\rho)-
\overline{T_k(\rho)}+T_k(n)-
\overline{T_k(n)}\|_{L^{\Gamma_{min}+1}(Q_T/Q_{T,k})}^{\frac{\Gamma_{min}+1}{2\Gamma_{min}}}\\
\le& C \|T_k(\rho)-
\overline{T_k(\rho)}\|_{L^{1}(Q_T/Q_{T,k})}^{\frac{\Gamma_{min}-1}{2\Gamma_{min}}}+C\|T_k(n)-
\overline{T_k(n)}\|_{L^{1}(Q_T/Q_{T,k})}^{\frac{\Gamma_{min}-1}{2\Gamma_{min}}}.\end{split}
\ee

Note that \be\label{divugotozero}
\lim\limits_{k\rightarrow\infty}\|\mathrm{div}u\|_{L^2(Q_{T,k})}=0
\ee since the Lebesgue measure of $Q_{T,k}$ converges to zero as
$k\rightarrow\infty$, due to \bex
\int_{Q_T}\Big(n^{\Gamma+\theta_1}+\rho^{\gamma+\theta_2}\Big)\,dx\,dt\le
C \eex given by (\ref{3-hierho}).

Therefore, to get (\ref{3-c2}), it suffices to prove \bex
\|T_k(\rho)- \overline{T_k(\rho)}\|_{L^{1}(Q_T/Q_{T,k})}+\|T_k(n)-
\overline{T_k(n)}\|_{L^{1}(Q_T/Q_{T,k})}\rightarrow 0 \eex as
$k\rightarrow\infty$, according to (\ref{3-a15}) and
(\ref{3-a16}).

Recalling that $T_k(z)=z$ if $z\le k$, we have
 \be\label{3-a6}
\begin{split}
&\|T_k(\rho)-
\overline{T_k(\rho)}\|_{L^{1}(Q_T/Q_{T,k})}+\|T_k(n)-
\overline{T_k(n)}\|_{L^{1}(Q_T/Q_{T,k})}\\=&\|\rho-
\overline{T_k(\rho)}\|_{L^{1}(Q_T/Q_{T,k})}+\|n-
\overline{T_k(n)}\|_{L^{1}(Q_T/Q_{T,k})}\\
\le&\liminf\limits_{k\rightarrow\infty}\|\rho_\delta-T_k(\rho_\delta)\|_{L^{1}(Q_T)}+\|n_\delta-
T_k(n_\delta)\|_{L^{1}(Q_T)}\\
=&\liminf\limits_{k\rightarrow\infty}\|\rho_\delta-T_k(\rho_\delta)\|_{L^{1}(Q_T\cap\{\rho_\delta\ge
k\})}+\|n_\delta- T_k(n_\delta)\|_{L^{1}(Q_T\cap \{n_\delta\ge
k\})}\\
\le&C\liminf\limits_{k\rightarrow\infty}\|\rho_\delta\|_{L^{1}(Q_T\cap\{\rho_\delta\ge
k\})}+C\|n_\delta\|_{L^{1}(Q_T\cap \{n_\delta\ge k\})}\rightarrow
0
\end{split}
\ee as $\delta\rightarrow 0$, due to (\ref{3-hierho}).

Therefore we complete the proof of the lemma.

\endpf

\bigskip

Now we are ready to prove (\ref{3-claim}). In fact, in view of
(\ref{3-last4}) and (\ref{right side lemma}), we have
\be\label{3-19}\begin{split}
\lim\limits_{k\rightarrow\infty}\int_\Omega
[\overline{L_k(\rho)}-L_k(\rho)+\overline{L_k(n)}-L_k(n)]\,dx
\le0.
\end{split}
\ee

By the definition of $L(\cdot)$, it is not difficult to justify
that \be\label{3-last6}\begin{cases}
\lim\limits_{k\rightarrow\infty}\Big[\|L_k(\rho)-\rho\log\rho\|_{L^1(\Omega)}+\|L_k(n)-n\log n\|_{L^1(\Omega)}\Big]=0,\\[2mm]
\lim\limits_{k\rightarrow\infty}\Big[\|\overline{L_k(\rho)}-\overline{\rho\log\rho}\|_{L^1(\Omega)}
+\|\overline{L_k(n)}-\overline{n\log n}\|_{L^1(\Omega)}\Big]=0.
\end{cases}
\ee

(\ref{3-19}) and (\ref{3-last6}) yields
\be\label{3-20}\begin{split} \int_\Omega
[\overline{\rho\log\rho}-\rho\log\rho+\overline{n\log n}-n\log
n]\,dx \le 0.
\end{split}
\ee On the other hand, since $\rho\log\rho\le
\overline{\rho\log\rho}$ and $n\log n\le\overline{n\log n}$ due to
the convexity of $z\mapsto z\log z$, we have \bex
\overline{\rho\log\rho}=\rho\log\rho \quad \mathrm{and }\;\;\,\
\overline{n\log n}=n\log n. \eex It allows us to have the
 strong
convergence of $\rho_\delta$ and $n_\delta$ in $L^{\gamma_1}(Q_T)$
and in $L^{\Gamma_1}(Q_T)$ for any
$\gamma_1\in[1,\gamma+\theta_2)$ and
$\Gamma_1\in[1,\Gamma+\theta_1)$, respectively. Therefore we
proved (\ref{3-claim}).

Then the proof of Theorem \ref{th:1.1} is complete.

\section*{Acknowledgements}
The main part of the work was done when the author visited the Department of
Energy and Petroleum Engineering, University of Stavanger, Norway. He would like to thank the
department for its kind hospitality. He also would like to
thank Professors Steinar Evje, Alexis Vasseur, Lei Yao, Cheng Yu, and Changjiang Zhu for some inspiring
discussions about the model. This work was partially supported by the National Natural Science
Foundation of China (Grant No. 12071152 and 11722104), and by the Guangdong Basic and Applied Basic Research Foundation (Grant No. 2020B1515310015).


\addcontentsline{toc}{section}{\\References}

\begin{thebibliography}{00}


\bibitem{BLS} J.-W. Barrett, Y. Lu, E. Suli.
Existence of large-data finite-energy global weak solutions to a
compressible Oldroyd-B model. Commun. Math. Sci., 15 (2017), 1265-1323.



\bibitem{Br} C. E. Brennen, Fundamentals of Multiphase Flow, Cambridge Univ. Press, 2005.

\bibitem{BDGG}D. Bresch, B. Desjardins, J.-M. Ghidaglia, E. Grenier.
Global weak solutions to a generic two-fluid model. Arch. Rational
Mech. Anal., 196 (2010), 599-629.

\bibitem{BHL} D. Bresch, X.D. Huang, J. Li. Global weak solutions to
one-dimensional nonconservative viscous compressible two-phase
system. Commun. Math. Phys., 309 (2012), 737-755.

\bibitem{BJ}
D. Bresch, P.-E. Jabin. Global existence of weak solutions for compressible Navier-Stokes equations: thermodynamically unstable pressure and anisotropic viscous stress tensor.
 Ann. of Math., (2) 188 (2018), no. 2, 577-684.

\bibitem{Bresch-Mucha}D. Bresch, P.-B. Mucha,
E. Zatorska. Finite-energy solutions for compressible two-fluid
Stokes system. Arch. Ration. Mech. Anal., 232 (2019), no. 2, 987-1029.


\bibitem{DL1} R. J. DiPerna, P.-L. Lions. Ordinary differential equations, transport theory and Sobolev spaces. Invent. Math., 98 (3) (1989), 511-547.

\bibitem{DL} R. J. DiPerna, P.-L. Lions. On the Cauchy problem for Boltzmann equations: global existence and weak stability. Ann. of Math., (2) 130 (1989), no. 2, 321-366.

\bibitem{E}S. Evje. Weak solution for a gas-iquid model relevant for describing gas-kick oil wells. SIAM J. Math. Anal., 43 (2011), 1887-1922.

\bibitem{EvjeCES2017}S. Evje. An integrative multiphase model for cancer cell migration under influence of
physical cues from the tumor microenvrionment. Chem Engineering Science, 204-259,
2017.

\bibitem{EK}S. Evje, K.H. Karlsen. Global existence of weak solutions for a viscous two-fluid model. J. Differential Equations, 245 (9) (2008), 2660-2703.

\bibitem{EK2} S. Evje, K.H. Karlsen. Global weak solutions for a viscous liquid-gas model with singular pressure law. Commun. Pure Appl. Anal., 8 (2009), 1867-1894.

\bibitem{EvjeSIAM}S. Evje, H.Y. Wen. A Stokes two-fluid model for cell migration that
can account for physical cues in the
microenvironment. SIAM J. Math. Anal.,
50 (2018), No. 1, 86-118.


\bibitem{EWZ} S. Evje, H.Y. Wen, C.J. Zhu. On global solutions to the viscous liquid-gas model with unconstrained transition to single-phase flow. Math. Models Methods Appl. Sci., 27 (2017), no. 2, 323-346.

\bibitem{Feiresil1}
E. Feireisl. On compactness of solutions to the compressible
isentropic Navier-Stokes equations when the density is not square
integrable. Comment. Math. Univ. Carolinae, 42 (1) (2001), 83-98.

\bibitem{FeireislJDE}E. Feireisl. Compressible Navier-Stokes equations with a non-monotone
pressure law. J. Differential Equations, 184 (2002), 97-108.

\bibitem{Feireisl2} E. Feireisl. Dynamics of Viscous
Compressible Fluids, Oxford University Press, Oxford, 2004.

\bibitem{Feireisl}
E. Feireisl, A. Novotny, H. Petzeltova. On the existence of
globally defined weak solutions to the Navier-Stokes equations. J.
Math. Fluid Mech., 3(2001), 358-392.

\bibitem{Feireisl-Novotny}E. Feireisl, A. Novotny. Singular Limits
in Thermodynamics of Viscous Fluids. Birkh\"auser Verlag, Basel,
2009.








\bibitem{Guo-Yang-Yao}Z. H. Guo, J. Yang, L. Yao. Global strong solution for a three-dimensional viscous
liquid-gas two-phase flow model with vacuum. J. Math. Phys., 52 (2011), 093102.


\bibitem{Hao-Li}C.C. Hao, H.L. Li. Well-posedness for a multidimensional viscous liquid-gas two-phase flow model. SIAM J.
Math. Anal., 44 (2012), 1304-1332.


\bibitem{Hu-Wang} X.P. Hu, D.H. Wang. Global existence and large-time behavior of solutions to the three-dimensional equations of compressible magnetohydrodynamic flows. Arch. Ration. Mech. Anal., 197 (2010), no. 1, 203-238.

\bibitem{I} M. Ishii. Thermo-Fluid Dynamic Theory of two-fluid Flow, Eyrolles, Paris, 1975.


\bibitem{Jiang-Zhang}S. Jiang, P. Zhang. Global spherically symmetry solutions of
the compressible isentropic Navier-Stokes equations. Comm. Math.
Phys., 215 (2001), 559-581.



\bibitem{Li-Sun}Y. Li, Y.Z. Sun. Global weak solutions to a two-dimensional compressible MHD equations of viscous non-resistive fluids. J. Differential Equations, 267(6)(2019), 3827-3851.

\bibitem{Lions}
P. L. Lions. Mathematical Topics in Fluid Mechanics, vol. II,
Compressible Models, Clarendon Press, Oxford, 1998.





\bibitem{MMMNPZ} D. Maltese, M. Michek, P. Mucha, A. Novotn\'y, M. Pokorn\'y, E. Zatorska. Existence of weak solutions for compressible Navier-Stokes equations with entropy transport. J. Differential Equations, 261 (2016), no. 8, 4448-4485.


\bibitem{MV3} A. Mellet, A. Vasseur. Asymptotic analysis for a Vlasov-Fokker-Planck/compressible Navier-Stokes system of equations. Comm. Math. Phys., 281 (2008), no. 3, 573-596.

\bibitem{Novotny}A. Novotn\'y, M. Pokorn\'y. Weak solutions for some compressible
multicomponent fluid models. Arch. Ration. Mech. Anal., 235(2020), 355-403.



\bibitem{Vasseur-Wen-Yu} A. Vasseur, H.Y. Wen, C. Yu. Global weak solution to the viscous two-fluid model with finite
energy. J. Math. Pures Appl., (9) 125 (2019), 247-282.


\bibitem{Yao-Zhang-Zhu} L. Yao, T. Zhang, C.J. Zhu. Existence and asymptotic behavior of global weak solutions to a 2D viscous liquid-gas
two-fluid flow model. SIAM J. Math. Anal., 42 (2010), 1874-1897.

\bibitem{Yao-Zhu}
L. Yao, C.J. Zhu. Free boundary value problem for a viscous
two-phase model with mass-dependent viscosity. J. Differential
Equations, 247 (2009), 2705-2739.

\bibitem{Yao-Zhu2}
L. Yao, C.J. Zhu. Existence and uniqueness of global weak solution
to a two-phase flow model with vacuum. Math. Ann., 349 (2011),
903-928.

\bibitem{Z}N. Zuber. On the dispersed two-phase flow in the laminar flow regime. Chem. Engrg. Sci., 19 (1964), 897-917.




\end{thebibliography}
\end{document}